\documentstyle[12pt]{article}

\newlength{\abstractwidth}
\renewcommand{\thefootnote}{\fnsymbol{footnote}}
\renewcommand{\thanks}[1]{\footnote{#1}} % Use this for footnotes
\newcommand{\starttext}{ \setcounter{footnote}{0}
\renewcommand{\thefootnote}{\arabic{footnote}}}

\newcommand{\be}{\begin{equation}}
\newcommand{\bea}{\begin{eqnarray}}
\newcommand{\eea}{\end{eqnarray}} 
\newcommand{\ee}{\end{equation}}

\def\ba{\begin{eqnarray}}
\def\ea{\end{eqnarray}}

\def\o{\omega}
\def\Re{{\rm Re}}

\def\log{\,{\rm log}\,}

\def\o{\omega}

\def\e{\varepsilon}

\def\o{\omega}

\def\p{\partial}

\def\Z{{\bf Z}}

\def\R{{\bf R}}
\def\C{{\bf C}}

\def\ddb{{\partial\bar\partial}}

\def\F{{\cal F}}

\def\[{{\bf [}}
\def\]{{\bf ]}}

\begin{document}
\starttext \baselineskip=15pt \setcounter{footnote}{0}
\newtheorem{theorem}{Theorem}
\newtheorem{lemma}{Lemma}
\newtheorem{definition}{Definition}
\newtheorem{proposition}{Proposition}
\newtheorem{corollary}{Corollary}

\begin{center}
{\Large \bf THE FU-YAU EQUATION WITH NEGATIVE SLOPE PARAMETER \footnote{Work supported in part by the National Science Foundation under Grant DMS-1266033, DMS-1605968 and DMS-1308136. Keywords: Hessian equations
with gradient terms; symmetric functions of eigenvalues; Moser iteration; maximum principles. 
%AMS classification numbers: 32Q26 (32Q15, 32Q20, 32U05, 32W20), 35Kxx.
}}

\bigskip

{\large Duong H. Phong, Sebastien Picard, and Xiangwen Zhang} \\

\medskip
\begin{abstract} 

{\small The Fu-Yau equation is an equation introduced by J. Fu and S.T. Yau as a generalization to arbitrary dimensions of an ansatz for the Strominger system. As in the Strominger system, it depends on a slope parameter $\alpha'$. The equation was solved in dimension $2$ by Fu and Yau in two successive papers for $\alpha'>0$, and for $\alpha'<0$. In the present paper, we solve the Fu-Yau equation 
in arbitrary dimension for $\alpha'<0$. To our knowledge, these are the first non-trivial solutions of the Fu-Yau equation in any dimension strictly greater than $2$.}

\end{abstract}
\end{center}

\section{Introduction}

\par

Let $(X,\o)$ be a compact K\"ahler manifold of dimension $n$. The Fu-Yau equation with slope parameter $\alpha'$ is the following equation for an unknown scalar function $u$,
\be \label{FY_eqn}
i \ddb (e^u \o - \alpha' e^{-u} \rho) \wedge \o^{n-2} + n \alpha' i \ddb u \wedge i \ddb u \wedge \o^{n-2} + \mu {\o^n \over n!} =0. 
\ee
Here $\mu : X \rightarrow {\bf R}$ is a smooth function satisfying $\int_X \mu {\o^n\over n!}= 0$, $\rho$ is a smooth real $(1,1)$ form, and the solution $u$ is required to be admissible, in the sense that the vector $\lambda'$
of eigenvalues of the Hermitian form
\be
\label{g_prime}
\o' = e^u \o + \alpha' e^{-u} \rho + 2n\alpha' i\partial\bar\partial u,
%g_{\bar kj}'=e^u g_{\bar kj}+\alpha' e^{-u}\rho_{\bar kj}+2n\alpha' u_{\bar kj}
\ee
with respect to $\o$, lies in the admissible cone $\Gamma_2$ defined in (\ref{defcone}) below. Henceforth, to simplify the notation, we shall just denote $\alpha'$ by $\alpha$. When $\alpha=0$, the Fu-Yau equation reduces to a Laplacian equation in $e^u$, so the only non-trivial cases are when $\alpha$ is strictly positive or negative. The Fu-Yau equation was solved in dimension ${\rm dim}\,X=2$ by Fu and Yau in two ground breaking papers, first for $\alpha>0$ in \cite{FY1}, and then for $\alpha<0$ in \cite{FY2}. The main goal of the present paper is to prove the following theorem:

\begin{theorem} \label{main_theorem}
Let $\alpha < 0$. Then for any dimension $n\geq 2$,
any smooth $(1,1)$-form $\rho$ and any smooth function $\mu$ satisfying the condition $\int_X\mu\,\o^n=0$, there exists a constant $M'$ so that, for all $M_0\geq M'$, there exists a smooth
admissible solution $u$ to the Fu-Yau equation (\ref{FY_eqn})
with normalization
\be
\int_X e^u=M_0.
\ee
\end{theorem}

\bigskip
The Fu-Yau equation in general dimension $n$ and with $\alpha>0$ has been studied in \cite{PPZ3}. All the basic a priori estimates needed for a solution by the method of continuity had been derived there, except for a lower bound on the second symmetric function of the eigenvalues of $\o'$. As a consequence, whether the equation is solvable for $\alpha>0$ and $n> 2$ is still an open question at this time. It is an intriguing question whether the distinct behavior of the equations with $\alpha>0$ and $\alpha<0$ is indicative of a significant geometric difference in the metrics defined by the two equations.

\medskip

The Fu-Yau equation is motivated by the fact that, in dimension $n=2$, as shown in \cite{FY1}, it is equivalent to the Strominger system for a certain class of $3$-dimensional manifolds constructed by Goldstein and Prokushkin \cite{GP}. In higher dimensions, it corresponds to an interesting modification of the Strominger system. More precisely, 
let $M$ be an $(n+1)$-dimensional complex manifold, equipped 
with a nowhere vanishing holomorphic $(n+1,0)$ form $\Omega$. Let $E \rightarrow M$ be a holomorphic vector bundle with Hermitian metric $H$. The {\it modified} Strominger system, as proposed by Fu and Yau \cite{FY1},
is the following system for a metric $\o$ on $M$
and a Hermitian metric $H$ on $E$,
\be \label{HE}
F_H \wedge \o^{n} = 0, \ \ F_H^{2,0} = F_H^{0,2} = 0
\ee
\be
\label{anomaly}
\bigg\{ i \ddb \o - {\alpha \over 4} ( {\rm Tr}\, R \wedge R - {\rm Tr}\, F_H \wedge F_H) \bigg\} \wedge \o^{n-1} = 0,
\ee
\be \label{balanced}
d \bigg( \| \Omega \|_{\o}^{{2(n-1) \over n}} \o^{n} \bigg) = 0.
\ee
Here $F_H$ is the curvature of the bundle $E$ with respect to the metric $H$ and $R$ is the Riemann curvature tensor of the Chern connection of the metric $\omega$, viewed as $(1,1)$-forms valued in the endomorphisms of $E$ and $T^{1,0}(M)$ respectively. Note that a natural extension to arbitrary dimensions of the Strominger system may have been with the power $\o^{n-2}$ in the equation (\ref{anomaly}). The power $\o^{n-1}$ above is motivated by the requirement that the modified Strominger system be equivalent to the Fu-Yau equation when restricted to Goldstein-Prokushkin fibrations\footnote{In \cite{FY1}, there appears to be a misprint, with the power of $\o$ in (\ref{anomaly}) written as $n-2$. The authors are grateful to Li-Sheng Tseng for pointing out to them the correct power $n-1$.}. When $n+1=3$, which is the case of particular interest in string theory, both powers $\o^{n-1}$ and $\o^{n-2}$ lead to the same Fu-Yau equation when restricted to Goldstein-Prokushkin fibrations.
As observed by Li and Yau \cite{LY}, in dimension $n+1=3$ the equation (\ref{balanced}) can be replaced by the equation $d^\dagger\o=i(\bar\partial-\partial)\ln \,\|\Omega\|_\o$, which is the original form written down by Strominger \cite{S}. Strominger's original motivation was from string theory, and the system he proposed would guarantee the $N=1$ supersymmetry of the heterotic string compactified to Minkowski space-time by a 3-dimensional internal space $M$. For this, we would take $n+1=3$, and the slope $\alpha$ is positive. Nevertheless, the Strominger systems and their modifications for general values of $n$ and $\alpha$ are compelling systems of considerable geometric interest, as they unify in a natural way two basic equations of complex geometry, namely the Hermitian-Einstein equation for holomorphic vector bundles and a generalization of the Ricci-flat equation. In particular, the equations (\ref{anomaly}) and (\ref{balanced}) (for fixed metric $H$ on $E$) can then be viewed as legitimate non-K\"ahler alternatives to the canonical metrics of K\"ahler geometry. 

\medskip
The way a solution of the Fu-Yau equation would give rise to a solution of the Strominger system is a key contribution of Fu and Yau \cite{FY1,FY2}, based on an earlier geometric construction of Goldstein and Prokushkin \cite{GP}.  Recall that the Goldstein-Prokushkin construction associates a toric fibration
$\pi: M\to X$ to a compact Calabi-Yau manifold $(X,\o_X)$ of dimension $n$ with nowhere vanishing holomorphic $n$-form $\Omega_X$ and two primitive harmonic forms ${\o_1\over 2\pi},{\o_2\over 2\pi}\in H^2(X,\Z)$. Furthermore, 
there is a $(1,0)$-form $\theta$ on $M$ so that $\Omega=\Omega_X\wedge\theta$ is a nowhere vanishing holomorphic $(n+1)$-form on $M$, and
\bea
\o_u=\pi^*(e^u\o_X)+i\theta\wedge\bar\theta
\eea
is a metric satisfying the balanced condition (\ref{balanced}) for any scalar function $u$ on $X$. Let $E\to X$ be a stable holomorphic vector bundle of degree $0$, and let $H$ be a Hermitian-Einstein metric on $E$, which exists by the Donaldson-Uhlenbeck-Yau theorem. Let $\pi^*(E)$, $\pi^*(H)$ be their pull-backs to $M$. The equations (\ref{HE}) and (\ref{balanced})
are now satisfied. It is then shown by Fu and Yau \cite{FY1, FY2} that the 
last equation (\ref{anomaly}) in the Strominger system for the system $(\pi^*E,\pi^*H, M,\o_u)$ is satisfied if and only if $u$ satisfies the Fu-Yau equation (\ref{FY_eqn}) on the manifold $X$, for suitable $\rho$ and $\mu$ given explicitly by
\be
\label{mu}
\mu {\o_X^n \over n!} = {(n-2)! \over 2} (\|\o_1 \|^2_{\o_X} + \| \o_2 \|^2_{\o_X}) {\o^n_X \over n!} + {\alpha \over 4} {\rm Tr}( F_H \wedge F_H - R_X \wedge R_X) \wedge \o^{n-2}_X.
\ee
The solvability condition $\int_X \mu\,{\o_X^n\over n!}=0$ of the Fu-Yau equation can be viewed as a cohomological condition on $X,E,\o_1,\o_2$ and the slope parameter $\alpha$.
Applying this construction of Fu and Yau, we obtain, as an immediate corollary of Theorem \ref{main_theorem}:

\begin{theorem}
\label{2}
Let $(X,\o_X)$ be an $n$-dimensional Calabi-Yau manifold, equipped with a nowhere vanishing holomorphic $n$-form $\Omega$ ($n\geq 2$). Let ${\o_1\over 2\pi},{\o_2\over 2\pi}\in H^{1,1}(X,{\bf Z})$ be primitive harmonic forms. Let $E\to X$ be a stable bundle with Hermitian-Einstein metric $E$. Assume that $\alpha<0$, and the cohomological condition $\int_X \mu\,{\o_X^n\over n!}=0$ is satisfied. 
Then the modified Strominger system (\ref{HE}), (\ref{anomaly}), (\ref{balanced}) admits a smooth solution of the form
$(\pi^*E, \pi^*H, M, \o_u)$.
\end{theorem}

Following Fu-Yau \cite{FY1,FY2}, we can construct many examples of data $(X,E,\o_1,\o_2,\alpha)$ satisfying the cohomological condition $\int_X \mu\,{\o_X^n\over n!}=0$ with $\alpha<0$ in higher dimension, and this is illustrated in \S 7. Thus
Theorem \ref{2} provides the first known solutions of the Strominger system in higher dimensions by the Fu-Yau ansatz. 
In section \S 7, we shall exhibit a specific example due to Fu and Yau \cite{FY2} with $\alpha=-2$ and $\mu=0$. 
Other, more geometric, constructions of solutions to the Strominger system have  also been provided in \cite{AG, Fe1, Fe2, FeY, FIUV1, FIUV2, FTY, Gr, UV}.

\medskip
Besides its occurrence in geometry and physics, the Fu-Yau equation (\ref{FY_eqn}) is also interesting from the point of view of the theory of fully non-linear elliptic partial differential equations. In dimension $n=2$, it is a complex Monge-Amp\`ere equation, and the natural ellipticity condition is that the form $\o'$ defined in (\ref{g_prime}) be positive-definite. But in dimension $n>2$, it is actually a complex $2$-Hessian equation, with the ellipticity condition given by the condition that the eigenvalues of $\o'$ be in the cone $\Gamma_2$. In fact, as worked out in detail in \S \ref{g'_section}, it can be rewritten as
\begin{eqnarray}\label{2hessian}
\frac{(\o')^2 \wedge \omega^{n-2}}{\o^n}=\frac{n(n-1)}{2}\left(e^{2u} - 4\alpha e^u |Du|^2 \right) + \nu,
\end{eqnarray}
where $\o'$ is the Hermitian $(1, 1)$ form given in (\ref{g_prime}), and
$\nu$ is a function depending on $u, Du$, $\mu$ and $\rho$, given explicitly in (\ref{kappa_equation}) below. Complex Hessian equations on compact manifolds have been studied extensively by many authors in recent years, see for example, \cite{Bl, DK, DK1, GS, HMW, KN, Lu, LN, Sun, Sz, STW, DZhang, XZhang}. However, in comparison with the previous works, a crucial new feature of the equation here is the dependence of the right hand side on the gradient $Du$ of the unknown function $u$.

\medskip
The proof of Theorem \ref{main_theorem} is by the method of continuity, and the main task is to derive the a priori estimates. We now describe briefly some of the innovations required in the derivation of these estimates. 

\smallskip
In the original papers of Fu-Yau \cite{FY1,FY2}, the $C^0$ estimate for equation (\ref{FY_eqn}) was proved using two different arguments depending on the sign of $\alpha$. Here we provide a unified approach. We also impose a simpler normalization condition on $\int_X e^u$ instead of on $\int_X e^{-pu}$, with $p$ depending on dimension $n$, as in \cite{FY1, FY2}.
This simpler normalization arises naturally in the study of a parabolic version of the Fu-Yau equation, which is a reduction by the Goldstein-Prokushkin construction of a  
geometric flow on $(2,2)$ forms introduced by the authors \cite{PPZ3} and called the anomaly flow. The anomaly flow preserves the balanced condition of metrics, and its stationary points satisfy the anomaly equation of the Strominger system. It can be shown that $\int_X e^{u}$ is constant along the flow.

\smallskip

The $C^1$ estimate uses the blow-up and Liouville theorem technique of Dinew-Kolodziej \cite{DK}. To adapt this argument, we need to show that {\it there is a uniform constant $C$, depending only on $\o, \rho$, $\mu$ and $\alpha$, such that}
\begin{eqnarray}\label{C2estimate}
\sup_X | \p \bar{\p} u |_{\o} \leq C(1 + \sup_X | Du |_{\o}^2).
\end{eqnarray}
For standard complex Hessian equations on compact K\"ahler manifolds, this $C^2$ estimate was obtained by Hou-Ma-Wu \cite{HMW}. It is worth mentioning that such a $C^2$ estimate combined with a blow-up argument has been a key ingredient in the solvability of several equations in complex geometry. For example, this type of estimate can be obtained for the form-type Monge-Amp\`ere equation occurring in the proof of the Gauduchon conjecture \cite{STW, TW}, and the dHYM equation \cite{CJY} motivated from mirror symmetry. In this paper, we establish the $C^2$ estimate (\ref{C2estimate}) for the Fu-Yau equation (\ref{FY_eqn}). The proof is quite different from the proof given by Fu and Yau \cite{FY1,FY2} for the case $n=2$.

\smallskip 

To obtain estimate (\ref{C2estimate}), the major obstacle in our case is the presence of gradient terms such as $e^u|Du|^2$. Indeed, if we allow the constant $C$ to depend on the gradient of $u$, the $C^2$ estimate for complex Hessian equations with gradient terms was established in \cite{PPZ2} under a stronger cone condition, building on the techniques developed by Guan-Ren-Wang \cite{GRW} for the real Hessian equations. However, the argument is not completely applicable here, because we need to get estimate (\ref{C2estimate}) with $C$ independent of the gradient in order to get subsequently the $C^1$ estimate by blow-up arguments. In this paper, we exploit the precise form of the gradient terms to obtain a cancellation of the $\sum |u_{ip}|^2$ terms arising from differentiating the right-hand side, as seen in Lemma \ref{c2lemma} in section \S \ref{C2estimatesec}. 

\smallskip

We would like to stress that the estimates here are also quite different from the case $\alpha>0$ which was studied in our previous work \cite{PPZ1}. When $\alpha>0$, the $C^1$ estimate easily follows from the $C^0$ estimate and ellipticity of the equation. However, a new difficulty arises, which is that {\it the equation may become degenerate}. This would happen if  the right-hand side $c_n (e^{2u} - 4 \alpha e^u |Du|^2)$ of the equation (\ref{2hessian}) tends to $0$. This is not possible if $\alpha$ is negative, but it cannot be ruled out at the outset if $\alpha$ is positive. In \cite{PPZ1}, we reduced the solvability of the Fu-Yau equation with $\alpha>0$ to a non-degeneracy estimate. In fact, this non-degeneracy estimate is equivalent to a strong Fu-Yau type gradient estimate, which was obtained in \cite{FY1} for $n=2$. It is still not known whether it holds for general $n$. In any case, a proof will certainly require a new method.

\medskip

The paper is organized as follows. We give a general setup for the method of continuity in section \S 2 and establish the a priori estimates in sections \S 3, \S 4, \S 5. In section \S 6, we give the proof of Theorem 1. In section \S 7, we follow Fu-Yau's construction to give a solution with $\alpha=-2$ to the modified Strominger system. Finally, in section \S 8, we propose another possible generalization of the Strominger system to higher dimensions which makes use instead of higher Chern classes.

\section{The Continuity Method}
\setcounter{equation}{0}

We shall work on a compact K\"ahler manifold $(X,\o)$ of dimension $n\geq 2$. We use the notation $\o=\sum g_{\bar kj} \, idz^j\wedge d\bar z^k$ and $\rho = \sum \rho_{\bar kj} \, idz^j\wedge d\bar z^k$, and we normalize the volume ${\rm Vol}(X,\o)=\int_X {\o^n\over n!}$ to be $1$. We denote by $D$ the covariant derivatives with respect to the background metric $\o$. All norms and inner products are with respect to the background metric $\o$ unless denoted otherwise.

\subsection{The set-up for the continuity method}

We shall solve the Fu-Yau equation using the continuity method. For a real parameter $t$, we consider
\be \label{continuity_meth1}
i \ddb (e^u \o - t \alpha e^{-u} \rho) \wedge \o^{n-2} + n \alpha i \ddb u \wedge i \ddb u \wedge \o^{n-2} + t \mu {\o^n \over n!} =0. 
\ee
Define $\lambda'_{(t,u)}$ to be the eigenvalues of 
\be
(g'_{(t,u)})_{\bar{j} k} = e^u g_{\bar{j} k} +t \alpha e^{-u} \rho_{\bar{j} k} + 2n \alpha u_{\bar{j} k}
\ee
with respect to the background K\"ahler metric $\o$. The tensor $(g'_{(t,u)})$ is relevant because the ellipticity condition of (\ref{continuity_meth1}) is $\lambda'_{(t,u)} \in \Gamma_2$, where 
\be
\label{defcone}
\Gamma_2=\{\lambda\in {\bf R}^n; \ \sigma_1(\lambda)>0,\ \sigma_2(\lambda)>0\}.
\ee
Here $\sigma_k(\lambda)$ is the $k$-th symmetric function of $\lambda$,
defined to be $\sigma_k(\lambda)=\sum_{j_1<\cdots<j_k}\lambda_{j_1}\cdots\lambda_{j_k}$. When $\lambda$ is the vector of eigenvalues of a Hermitian form $\o'$ with respect to the background form $\o$, $\sigma_k$ can also be expressed directly in terms of $\o'$ by
\be
\sigma_k(\o')={n\choose k} 
{(\omega')^k\wedge \omega^{n-k}\over\omega^n}.
\ee
It is a well-known property of the cone $\Gamma_2$ that for each index $l\in \{1, \cdots, n\}$, we have $\sum_{i \neq l} \lambda'_i >0$.

\smallskip 
Let $M_0 >0$ be a constant which will eventually be taken to be very large. For $0<\gamma<1$, we define the following function spaces
\be \label{cont_meth_normalization}
B_M = \{ u \in C^{2,\gamma}(X,\R) : \int_X e^{u} = M_0 \},
\ee
\be
B_1 = \{ (t, u) \in [0,1] \times B_M : \lambda'_{(t,u)} \in \Gamma_2 \},
\ee
\be
B_2 = \{ \psi \in C^\gamma(X,\R) : \int_X \psi = 0 \}.
\ee
Consider the operator $\Psi: B_1 \rightarrow B_2$ defined by
\be
\Psi(t,u) {\o^n \over n!} = i \ddb (e^u \o - t \alpha e^{-u} \rho) \wedge \o^{n-2} + n \alpha i \ddb u \wedge i \ddb u \wedge \o^{n-2} + t \mu {\o^n \over n!}.
\ee
Define the set 
\be
I = \{ t \in [0,1] : \ {\rm there \ exists} \ u \in B_M \ {\rm such \ that} \ (t, u) \in B_1 \ {\rm and} \ \Psi(t,u)=0 \}.
\ee
We note that $0 \in I$ due to the trivial solution $u_0 = \log M_0$. The goal is to show that $I$ is both open and closed. In the remaining part of this section, we show that $I$ is open. The proof of closedness and the necessary a priori estimates will be given in subsequent sections.

\subsection{Proof of the openness of $I$}

Suppose $\hat{t} \in I$. Then there exists $\hat{u} \in B_M$ such that $\Psi(\hat{t}, \hat{u})=0$. We wish to use the implicit function theorem to solve $\Psi(t,u_t)$ for $t$ close to $\hat{t}$. We compute the linearized operator at $\hat{u}$ to be
\be
(D_u \Psi)_{(\hat{t},\hat{u})} : \bigg\{ h \in C^{2,\gamma}(X,\R) : \int_X h e^{\hat{u}} = 0 \bigg\} \rightarrow \bigg\{\psi \in C^\gamma(X,\R) : \int_X \psi =0 \bigg\}, \ \ (D_u \Psi)_{(\hat{t},\hat{u})} = L,
\ee
with
\be
L(h) {\o^n \over n!} = i \ddb ( h e^{\hat{u}} \o + \hat{t} \alpha h e^{-\hat{u}} \rho) \wedge \o^{n-2} + 2 n \alpha i \ddb \hat{u} \wedge i \ddb h \wedge \o^{n-2}.
\ee
Expanding terms gives
\bea
L(h)  &=&  (n-2)! \, g^{i \bar{k}} g^{p \bar{j}} \tilde{g}_{\bar{k} p} D_i D_{\bar{j}} h + 2\Re \bigg\{ {i \p h \wedge \bar{\p} (e^{\hat{u}} \o + \hat{t} \alpha e^{-\hat{u}} \rho) \wedge \o^{n-2} \over (n!)^{-1} \o^n} \bigg\}  \nonumber\\
&& + \bigg\{ {i \ddb (e^{\hat{u}} \o + \hat{t} \alpha e^{-\hat{u}} \rho) \wedge \o^{n-2} \over  (n!)^{-1} \o^n} \bigg\} h,
\eea
where $\tilde{g}_{\bar{k} p} = (g^{a \bar{b}} (g'_{(\hat{t},\hat{u})})_{\bar{b} a}) g_{\bar{k} p} - (g'_{(\hat{t},\hat{u})})_{\bar{k} p} >0$ since $\lambda'_{(\hat{t},\hat{u})} \in \Gamma_2$. We see that $L$ is a second order linear elliptic operator on $X$, and the deformation $tL + (1-t)\Delta$ shows that $L$ has the same index as the Laplacian, namely index $0$. We shall compute the adjoint $L^*$ with respect to the $L^2$ inner product with volume ${\o^n \over n!}$. We integrate by parts to obtain
\bea
\int_X \psi \, L(h) \, {\o^n \over n!}  &=& \int_X \psi i \ddb ( h e^{\hat{u}} \o + \hat{t} \alpha h e^{-\hat{u}} \rho ) \wedge \o^{n-2} +  2 n \alpha \int_X \psi i \ddb \hat{u} \wedge i \ddb h \wedge \o^{n-2} \nonumber\\
&=& \int_X h \, \{e^{\hat{u}} \o + \hat{t} \alpha e^{-\hat{u}} \rho  +  2 n \alpha  i \ddb \hat{u} \} \wedge i \ddb \psi \wedge \o^{n-2} \nonumber\\
&=& (n-2)! \int_X h \, g^{i \bar{k}} g^{p \bar{j}} \tilde{g}_{\bar{k} p} D_i D_{\bar{j}} \psi \, {\o^n \over n!}.
\eea
It follows that $L^* = (n-2)! \, g^{i \bar{k}} g^{p \bar{j}} \tilde{g}_{\bar{k} p} D_i D_{\bar{j}}$. By the strong maximum principle, the kernel of $L^*$ is the constant functions. Again by the strong maximum principle, a non-zero function in the image of $L^*$ must change sign. Since $L^*$ has index 0, the codimension of ${\rm Im} \, L^*$ is one, and so the kernel of $L$ is spanned by a function of constant sign. To summarize, we have
\be
{\rm Ker} \,L^* = \R, \ \ {\rm Ker}\, L = \R \langle \phi \rangle, \ \ \phi \ {\rm constant \ sign}.
\ee
By the Fredholm alternative, we obtain that $(D_u \Psi)_{(\hat{t}, \hat{u})}$ is an isomorphism of tangent spaces. By the implicit function theorem, we can solve $\Psi(t,u_t)$ for $t$ close to $\hat{t}$. Hence $I$ is open.

\section{The $C^0$ Estimate}
\setcounter{equation}{0}

This section as well as the next three are devoted to the proof of the closedness of the set $I$ of parameters $t$ for which the deformed equation
(\ref{continuity_meth1}) can be solved. For this, we need a priori $C^0$, $C^1$, $C^2$, and $C^{2,\gamma}$ a priori estimates. As usual, it is notationally convenient to derive these bounds for the original equation (\ref{FY_eqn}), as long as the bounds obtained depend only on suitable norms for the data $\rho$ and $\mu$.

\subsection{The supremum estimate}

\begin{proposition} \label{sup_estimate}
Let $u$ be a solution to (\ref{FY_eqn}) such that $\lambda' \in \Gamma_2$ and $\int_X e^u = M_0$. Suppose $e^u \geq 1$ and $\alpha e^{-2u} \rho \geq - {1 \over 2} \o$. Then there exists a constant $C$ depending only on $(X,\o)$, $\alpha$, $\rho$ and $\mu$ such that
\be
\sup_X e^u \leq C \int_X e^u = C M_0.
\ee
\end{proposition}
{\it Proof.} We proceed by Moser iteration. Recall the form $\o'$ defined by
\be
\o' = e^u \o + \alpha  e^{-u} \rho  + 2n \alpha i \partial \bar{\partial} u.
\ee
The starting point is to compute the quantity
\be
\int_X i \ddb (e^{-ku}) \wedge \o' \wedge \o^{n-2}
\ee
in two different ways. On one hand, by the definition of $\o'$ and Stokes' theorem, we have
\be
\int_X i \ddb (e^{-ku}) \wedge \o' \wedge \o^{n-2} = \int_X \{ e^u \o + \alpha e^{-u} \rho \} \wedge i \ddb (e^{-ku}) \wedge \o^{n-2}.
\ee
Expanding
\bea \label{int_identity1}
\int_X i \ddb (e^{-ku}) \wedge \o' \wedge \o^{n-2} &=& k^2 \int_X e^{-ku} \{ e^u \o  + \alpha e^{-u} \rho \} \wedge i \p u \wedge \bar{\p} u \wedge \o^{n-2} \nonumber\\
&& - k \int_X e^{-ku} \{ e^u \o  + \alpha e^{-u} \rho \} \wedge i \ddb u \wedge \o^{n-2} . 
\eea
On the other hand, without using Stokes' theorem, we obtain
\bea \label{int_identity2}
\int_X i \ddb (e^{-ku}) \wedge \o' \wedge \o^{n-2}  &=& k^2 \int_X e^{-ku} i \p u \wedge \bar{\p} u \wedge \o' \wedge \o^{n-2} \nonumber\\
&&  -k \int_X  e^{-ku} \{ e^u\o + \alpha e^{-u} \rho \} \wedge i \ddb u \wedge \o^{n-2} \nonumber\\
&&- (2 n \alpha)k \int_X  e^{-ku} i \ddb u \wedge i \ddb u \wedge \o^{n-2} .
\eea
We equate (\ref{int_identity1}) and (\ref{int_identity2})
\bea
0 &=& - k^2 \int_X e^{-ku} i \p u \wedge \bar{\p} u \wedge \o' \wedge \o^{n-2} +k^2 \int_X e^{-ku} \{ e^u \o  + \alpha e^{-u} \rho \} \wedge i \p u \wedge \bar{\p} u \wedge \o^{n-2} \nonumber\\
&&+ (2 n \alpha)k \int_X  e^{-ku} i \ddb u \wedge i \ddb u \wedge \o^{n-2} .
\eea
Using equation (\ref{FY_eqn}),
\bea
0 &=&   - k^2 \int_X e^{-ku} i \p u \wedge \bar{\p} u \wedge \o' \wedge \o^{n-2} +k^2 \int_X e^{-ku} \{ e^u \o  + \alpha e^{-u} \rho \} \wedge i \p u \wedge \bar{\p} u \wedge \o^{n-2} \nonumber\\
&&- 2k \int_X e^{-ku} \mu {\o^n \over n!} -2k \int_X  e^{-ku} i \ddb (e^u \o - \alpha e^{-u} \rho) \wedge \o^{n-2}.
\eea
Expanding out terms and dividing by $2k$ yields
\bea
0 &=&   - {k \over 2} \int_X e^{-ku} i \p u \wedge \bar{\p} u \wedge \o' \wedge \o^{n-2} + {k \over 2} \int_X e^{-ku} \{ e^u \o  + \alpha e^{-u} \rho \} \wedge i \p u \wedge \bar{\p} u \wedge \o^{n-2} \nonumber\\
&&-  \int_X e^{-ku} \mu {\o^n \over n!} - \int_X e^{-(k-1)u} i \ddb u \wedge \o^{n-1} \nonumber\\
&&  - \int_X e^{-(k-1)u} i \p u \wedge \bar{\p} u \wedge \o^{n-1} - \alpha \int_X  e^{-(k+1)u} i \ddb u \wedge \rho \wedge \o^{n-2} \nonumber\\
&&+  \alpha \int_X  e^{-(k+1)u} i \p u \wedge \bar{\p} u \wedge \rho \wedge \o^{n-2} +  \alpha \int_X e^{-(k+1)u} i \ddb \rho \wedge \o^{n-2} \nonumber\\
&&- 2 \alpha \Re \int_X e^{-(k+1)u} i \p u \wedge \bar{\p} \rho \wedge \o^{n-2} .\nonumber
\eea
Integration by parts gives
\bea \label{integrals_byparts}
0 &=&    - {k \over 2} \int_X e^{-ku} i \p u \wedge \bar{\p} u \wedge \o' \wedge \o^{n-2} - {k \over 2} \int_X e^{-ku} \{ e^u \o  + \alpha e^{-u} \rho \} \wedge i \p u \wedge \bar{\p} u \wedge \o^{n-2} \nonumber\\
&& - \int_X e^{-ku} \mu {\o^n \over n!} +  \alpha \int_X e^{-(k+1)u} i \ddb \rho \wedge \o^{n-2} -  \alpha \int_X e^{-(k+1)u} i \p u \wedge \bar{\p} \rho \wedge \o^{n-2} .
\eea
One more integration by parts yields the following identity:
\bea \label{integrals_simplified}
&& {k \over 2} \int_X e^{-ku} \{ e^u \o + \alpha e^{-u} \rho \} \wedge i \p u \wedge \bar{\p} u \wedge \o^{n-2} \\&=& - {k \over 2}  \int_X e^{-ku} i \p u \wedge \bar{\p} u \wedge \o' \wedge \o^{n-2} - \int_X e^{-ku} \mu   + (\alpha - {\alpha \over k+1}) \int_X  e^{-(k+1)u}  i \ddb \rho \wedge \o^{n-2} . \nonumber 
\eea
The identity (\ref{integrals_simplified}) will be useful later to control the infimum of $u$, but to control the supremum of $u$, we replace $k$ with $-k$ in (\ref{integrals_simplified}). Then, for $k \neq 1$,
\bea \label{sup_key_identity}
&& {k \over 2} \int_X e^{(k+1)u} \{ \o + \alpha e^{-2u} \rho \} \wedge i \p u \wedge \bar{\p} u \wedge \o^{n-2} \\&=& - {k \over 2}  \int_X e^{ku} i \p u \wedge \bar{\p} u \wedge \o' \wedge \o^{n-2} + \int_X e^{ku} \mu   - (\alpha - {\alpha \over 1-k}) \int_X  e^{(k-1)u}  i \ddb \rho \wedge \o^{n-2} . \nonumber 
\eea
Since $\lambda' \in \Gamma_2$, by the properties of the cone we have that $\sum_{i \neq l} \lambda'_i >0$ for each index $l\in \{1, \cdots, n\}$. It follows that
\be \label{pos_estimate}
i \p u \wedge \bar{\p} u \wedge \o' \wedge \o^{n-2} \geq 0.
\ee
Let $\beta = {n \over n-1}>1$. We can use (\ref{pos_estimate}) and (\ref{sup_key_identity}) to derive the following estimate for any $k \geq \beta$
\be
k  \int_X e^{(k+1)u} \{ \o + \alpha e^{-2u} \rho \} \wedge i \p u \wedge \bar{\p} u \wedge \o^{n-2} \leq C \left( \int_X e^{ku} + \int_X e^{(k-1)u} \right).
\ee
By assumption, $\alpha e^{-2u} \rho \geq  - {1 \over 2} \o$. For $k \geq 2 \beta$, we can estimate
\be
 \int_X | D e^{{k \over 2}u}|^2 \leq C k \left( \int_X e^{(k-1)u} + \int_X e^{(k-2)u} \right).
\ee
Since $e^u \geq 1$, we can conclude
\be \label{sup_iter}
 \int_X | D e^{{k \over 2}u}|^2 \leq C k \int_X e^{ku},
\ee
for $k \geq 2 \beta$. The Sobolev inequality yields
\be \label{sup_iter}
\left( \int_X e^{k \beta u} \right)^{1/\beta} \leq C k \int_X e^{ku}.
\ee
After iterating this estimate, we arrive at
\be 
\sup_X e^{u} \leq C \|e^u\|_{L^{2 \beta}}.
\ee
To relate the $L^{2 \beta}$ norm of $e^u$ to $\int_X e^u = M_0$, we can use a standard scaling argument.
\be
\sup_{X} e^{u} \leq C \bigg( \int_{X} e^{u} e^{2\beta u -u} \bigg)^{1/2\beta} \leq   C \bigg( \sup_{X} e^{(2\beta -1)u} \bigg)^{1 \over 2 \beta} \bigg( \int_{X} e^{u} \bigg)^{1/2\beta}.
\ee
It follows immediately that
\be
\sup_X e^{u} \leq C \int_X e^u = CM_0.
\ee

\subsection{An integral estimate for $e^{-u}$}

\begin{proposition} \label{int_e^{-u}}
Let $u$ be a solution to (\ref{FY_eqn}) such that $\lambda' \in \Gamma_2$ and $\int_X e^u = M_0$. There exists $0<\delta'<1$, chosen small enough such that $e^{-u} \leq \delta'$ implies $\alpha e^{-2u} \rho \geq -{1 \over 2} \omega$, and $C_1$ depending only on $(X,\o)$, $\alpha$, $\rho$, $\mu$ with the following property. If $e^{-u} \leq \delta'$, then
\be
\int_X e^{-u} \leq {C_1 \over M_0}.
\ee
\end{proposition}
%\subsubsection{The Reverse Poincar\'e Inequality}
{\it Proof.} Setting $k=2$ in (\ref{integrals_simplified}) and using (\ref{pos_estimate}) gives 
\be
\int_X e^{-u} \{ \o + \alpha e^{-2u} \rho \} \wedge i \p u \wedge \bar{\p} u \wedge \o^{n-2}  \leq C \bigg( \int_X e^{-2u} + \int_X e^{-3u} \bigg) .
\ee
Choose $\delta'>0$ such that $\alpha e^{-2u} \rho \geq  -{1 \over 2} \o$. Since we are assuming that $e^{-u} \leq \delta'$ pointwise, we obtain
\be
\int_X |D e^{- {u \over 2}} |^2 \leq C \delta' \int_X e^{-u}.
\ee
By the Poincar\'e inequality
\be
\int_X e^{-u} - \bigg( \int_X e^{- {u \over 2}} \bigg)^2 = \int_X \bigg| e^{- {u \over 2}} - \int_X e^{- {u \over 2}} \bigg|^2 \leq C \int_X |D e^{- {u \over 2}} |^2.
\ee
Hence, for some $C_0$ independent of $\delta'$, if $\delta'$ is small
\be \label{reverse_poincare}
\int_X e^{-u} \leq {1 \over 1 - C_0 \delta'} \bigg( \int_X e^{- {u \over 2}} \bigg)^2.
\ee
%\subsubsection{The Measure Estimate}
Let $U = \{ e^{-u} \leq {2 \over M_0} \}$. Then by Proposition \ref{sup_estimate},
\be
M_0 = \int_U e^u + \int_{X \backslash U} e^u \leq |U| \sup_X e^u + (1 - |U|) {M_0 \over 2} \leq CM_0 |U| + (1- |U|) {M_0 \over 2}.
\ee
It follows that there exists $\theta>0$ independent of $M_0$ such that
\be
|U| > \theta >0 .
\ee
%\subsubsection{The Integral Estimate}
Let $\e>0$. We may use the measure estimate and (\ref{reverse_poincare}) to obtain
\bea
\bigg( \int_X e^{- {u \over 2}} \bigg)^2 &\leq& (1 + C_\e) \bigg( \int_U e^{- {u \over 2}} \bigg)^2 + (1 + \e) \bigg( \int_{X \backslash U} e^{- {u \over 2}} \bigg)^2 \nonumber\\
&\leq& (1+C_\e) |U| \int_U e^{-u} + (1+\e)(1-|U|) \int_{X \backslash U} e^{-u} \nonumber\\
&\leq& (1+C_\e) {2 \over M_0} + (1+\e) (1-\theta) {1 \over 1-C_0 \delta'} \bigg( \int_X e^{- {u \over 2}} \bigg)^2.
\eea
Thus
\be
\bigg( \int_X e^{- {u \over 2}} \bigg)^2 \leq (1+C_\e) {2 \over M_0} \bigg( {1 \over 1 - (1+\e) (1-\theta) (1-C_0 \delta')^{-1}} \bigg).
\ee
Therefore by (\ref{reverse_poincare})
\be
\int_X e^{-u} \leq {1 \over 1 - C_0\delta'} (1+C_\e) {2 \over M_0}  \bigg({1 \over 1 - (1+\e) (1-\theta)(1-C_0\delta')^{-1}}\bigg).
\ee
Choose $\e = {\theta \over 2}$ and suppose $0<\delta'<{\theta \over 4C_0}$. Using this choice of $\e$ and $\delta'$ along with $0<\theta<1$, it follows that $(1+\e)(1-C_0\delta')^{-1} \leq 1+\theta$. Therefore
\be
\int_X e^{-u} \leq {1 \over 1 - {\theta \over 4}} (1+C_\theta) {2 \over M_0} {1 \over\theta^2} = {C_1 \over M_0}.
\ee
The important point is that $\theta$ does not depend on $M_0$, so both $\delta'$ and $C_1$ do not depend on $M_0$.

\subsection{The infimum estimate}

\begin{proposition} \label{inf_estimate}
Let $u$ be a solution to (\ref{FY_eqn}) such that $\lambda' \in \Gamma_2$ and $\int_X e^u = M_0$. There exists $0<\delta'<1$ (the same $\delta'$ as in Proposition \ref{int_e^{-u}}) and $C_2$ depending only on $(X,\o)$, $\alpha$, $\rho$, $\mu$, such that if $e^{-u} \leq \delta'$, then
\be
\sup_X e^{-u} \leq {C_2 \over M_0}.
\ee
\end{proposition}
{\it Proof.} Combining (\ref{integrals_simplified}) with (\ref{pos_estimate}), and choosing $\delta'>0$ such that $\alpha e^{-2u} \rho \geq  -{1 \over 2} \o$, we obtain for $k \geq 2$
\be
k \int_X e^{-(k-1)u} |Du|^2 \leq C \bigg( \int_X e^{-ku} + \int_X e^{-(k+1)u} \bigg).
\ee
Therefore, for $k \geq 1$, we have
\be
\int_X | D e^{-{k \over 2}u}|^2 \leq C k \left( \int_X e^{-(k+1)u} + \int_X  e^{-(k+2)u} \right) . 
\ee
If $e^{-u}  \leq  \delta' <1$, we deduce
\be
\int_X | D e^{-{k \over 2}u}|^2 \leq C k \int_X e^{-ku}.
\ee
By the Sobolev inequality
\be
\left( \int_X e^{-k \beta u} \right)^{1/\beta} \leq C k \int_X e^{-ku}.
\ee
By iterating this estimate, we obtain
\be 
\sup_X e^{-u} \leq C \|e^{-u}\|_{L^{1}}.
\ee
Combining this estimate with Proposition \ref{int_e^{-u}} completes the proof.

\subsection{The $C^0$ estimate along the continuity method}

Combining the supremum and infimum estimates, we shall prove the desired $C^0$ estimate along the continuity method (\ref{continuity_meth1}).
\begin{proposition} \label{C0_est} Let $\alpha \neq 0$. There exists $B_1>1$, $B_2>1$, and $M' \gg 1$ depending only on $(X,\o)$, $\alpha$, $\rho$, and $\mu$ such that every $M_0 \geq M'$ has the following property. Let $u_0 = \log M_0$. Suppose that for all $t \in [0,t_0)$ with $t_0\leq 1$ there exists a solution $u_t$ to
\be
i \ddb (e^{u_t} \o - t \alpha e^{-{u_t}} \rho) \wedge \o^{n-2} + n \alpha i \ddb u_t \wedge i \ddb u_t \wedge \o^{n-2} + t \mu {\o^n \over n!} =0,
\ee
such that $\lambda'_{(t,u_t)} \in \Gamma_2$ and $\int_X e^{u_t}=M_0$. Then the following $C^0$ estimate holds:
\be
e^{u_t} \leq B_1 M_0, \ \ e^{-u_t} \leq {B_2 \over M_0}.
\ee
\end{proposition}
{\it Proof:} Choose $B_2 = 2 C_2$ where $C_2>1$ is as in Proposition \ref{inf_estimate}. Take $M' \gg 1$ such that $(2 C_2) M'^{-1} < \delta'$, where $\delta'>0$ is as in Proposition \ref{inf_estimate}. At $t=0$, we have $e^{-u_0} = M_0^{-1} < (2 C_2) M_0^{-1}$. We claim that $e^{-u_t}$ can never reach $(2 C_2) M_0^{-1}$ on $[0,t_0)$. If for some $t' \in (0,t_0)$ there holds $e^{-u_{t'}}= (2 C_2) M_0^{-1} < \delta'$, then by Proposition \ref{inf_estimate} it would follow that $e^{-u_{t'}} \leq C_2 M_0^{-1}$, which is a contradiction. The estimate on $e^{u_t}$ was established in Proposition \ref{sup_estimate}. Q.E.D.

\section{The $C^2$ estimate}
\setcounter{equation}{0}

We come now to one of the key estimates, namely the $C^2$ estimate. For this, it is essential to view the Fu-Yau equation as a complex $2$-Hessian equation to exploit the concavity of the operator.

\subsection{The Fu-Yau equation as a Hessian equation}\label{g'_section}

Using the elementary symmetric function, equation (\ref{FY_eqn}) can be written as the following scalar equation
\bea \label{scalar_FY_eqn}
0 &=& \{ (n-1)  e^u g^{j \bar{k}}  + \alpha e^{-u} \tilde{\rho}^{j \bar{k}} \} D_j D_{\bar{k}}u + 2n \alpha \sigma_2(i \ddb u) + (n-1) e^u |Du|^2 \nonumber\\
&&- \alpha e^{-u} \tilde{\rho}^{j \bar{k}} u_j u_{\bar{k}} - 2 \alpha \Re \langle \p e^{-u}, \p \rho \rangle_\o - \alpha e^{-u} \Delta_\o \rho  + {\mu \over (n-2)!}.
\eea
Here we introduced the following notation
\be
\tilde{\rho}^{j \bar{k}} = g^{j \bar{\ell}} g^{m \bar{k}} ((g^{a \bar{b}} \rho_{\bar{b} a} ) g_{\bar{\ell} m} - \rho_{\bar{\ell} m}),
\ee
\be
\langle \p e^{-u}, \p \rho \rangle_\o \, {\o^n \over n!} ={ i \p e^{-u} \wedge  \bar{\p} \rho \wedge \o^{n-2}\over  (n-2)! }, \ \ \Delta_\o \rho \, {\o^n \over n!} =  {i \ddb \rho \wedge \o^{n-2}\over  (n-2)! }.
\ee
As it was mentioned in (\ref{2hessian}), we shall rewrite the equation in terms of
\be
g'_{\bar{k} j}= e^u g_{\bar{k} j} + \alpha e^{-u} \rho_{\bar{k} j} +2n\alpha u_{\bar{k} j},
\ee
where $\o'=i\sum g'_{\bar k j} dz^j \wedge d \bar z^k$. We will use $\lambda'$ to denote the eigenvalues of $g'$ with respect to the background metric $\o$. Direct computation gives
\bea \label{sigma_2_lambda'}
\sigma_2(\lambda') &=& {n(n-1) \over 2} e^{2u}  + \alpha^2 e^{-2u} \sigma_2(\rho)+ (2n\alpha)^2 \sigma_2(i \ddb u) +  \alpha (n-1) g^{a \bar{b}} \rho_{\bar{b} a} \nonumber\\
&& + (2n \alpha) \{ (n-1)  e^u g^{j \bar{k}}  + \alpha e^{-u} \tilde{\rho}^{j \bar{k}} \} D_j D_{\bar{k}}u.
\eea
Introduce the constant $\kappa_c = {n(n-1) \over 2}$. When $\alpha \neq 0$, we may combine (\ref{scalar_FY_eqn}) and (\ref{sigma_2_lambda'}) to obtain the following equivalent equation
\bea \label{kappa_equation}
\sigma_2(\lambda') &=& \kappa_c (e^{2u} - 4 \alpha e^{u} |Du|^2) + 2n \alpha^2 e^{-u} \tilde{\rho}^{j \bar{k}} u_j u_{\bar{k}} - 4n \alpha^2 e^{-u} \Re \langle \p u, \bar{\p} \rho \rangle_\o  \nonumber\\
&& + \alpha^2 e^{-2u} \sigma_2 (\rho) + 2n\alpha^2 e^{-u} \Delta_\o \rho + \alpha (n-1) g^{a \bar{b}} \rho_{\bar{b} a} - {2n \alpha \over (n-2)!} \mu.
\eea
As noted in the introduction, the ellipticity condition is that $\lambda' \in \Gamma_2$.

\smallskip
To obtain higher order estimates, we will work with a version of the equation (\ref{kappa_equation}) which exhibits a concave elliptic operator. Denote $F = \sigma_2^{1/2}(\lambda')$. Equation (\ref{kappa_equation}) is equivalent to
\be
F=\sigma_2(\lambda')^{1/2} = w,
\ee
with
\bea \label{w_equation}
w^2 &=& \kappa_c e^{2u} -2 \alpha e^{u} \bigg\{2 \kappa_c |Du|^2 - n \alpha e^{-2u} \tilde{\rho}^{j \bar{k}} u_j u_{\bar{k}} + 2n \alpha e^{-2u} \Re \langle \p u, \bar{\p} \rho \rangle_\o \bigg\}  \nonumber\\
&& + \alpha^2 e^{-2u} \sigma_2 (\rho) + 2n\alpha^2 e^{-u} \Delta_\o \rho +\alpha(n-1)  g^{a \bar{b}} \rho_{\bar{b} a} - {2n \alpha \over (n-2)!} \mu.
\eea

\subsection{The linearization $F^{j \bar{k}}$}
At a point $p \in X$ where the background metric $g_{\bar{k} j} = \delta_{kj}$, we will use the notation 
\be
\sigma_2^{j \bar{k}} = {\p \sigma_2(\lambda') \over \p g'_{\bar{k} j}}, \ \  F^{j \bar{k}} = {\p \sigma_2^{1/2}(\lambda') \over \partial g'_{\bar{k} j}}, \ \ F^{j \bar{k}, \ell \bar{m}} = {\p^2 \sigma_2^{1/2}(\lambda') \over \partial g'_{\bar{k} j} \partial g'_{\bar{m} \ell}}.
\ee
Thus
\be
F^{j \bar{k}} = {\sigma_2^{j \bar{k}} \over 2 \sigma_2(\lambda')^{1/2}}, \ \ \F =  F^{j \bar{k}} g_{\bar{k} j} = {(n-1) \over 2} {\sigma_1(\lambda') \over  \sigma_2(\lambda')^{1/2}}.
\ee
In this section, we shall derive expressions for $F^{j \bar{k}} D_j D_{\bar{k}}$ acting on various quantities. First,
\bea \label{F_DDu}
2n\alpha F^{j\bar k}D_jD_{\bar k}u
&=&
F^{j\bar k}g_{\bar kj}'-e^u F^{j\bar k} g_{\bar kj} - \alpha e^{-u} F^{j\bar k} \rho_{\bar kj}
\nonumber\\
&=&
\sigma_2^{1/2}(\lambda')- e^u \F - \alpha e^{-u} F^{j\bar k} \rho_{\bar kj}.
\eea
Covariantly differentiating $\sigma_2(\lambda')^{1/2}$ gives 
\bea \label{D_sigma2}
\p_p \sigma_2^{1/2} =F^{j\bar k}D_pg_{\bar kj}'. 
\eea
Substituting in the definition of $g_{\bar kj}'$, we obtain the following formulas for $2n\alpha F^{j\bar k}D_jD_{\bar k}$ acting on $Du$,
\bea \label{LDu1}
2n\alpha F^{j\bar k}D_jD_{\bar k}(D_pu)=
\p_p \sigma_2^{1/2}
- F^{j \bar{k}} D_p ( e^u g_{\bar{k} j} + \alpha e^{-u} \rho_{\bar{k} j}),
\eea
\bea \label{LDu2}
2n\alpha F^{j\bar k}D_jD_{\bar k}(D_{\bar p}u)
=\p_{\bar p} \sigma_2^{1/2}
- F^{j \bar{k}} D_{\bar{p}} ( e^u g_{\bar{k} j} + \alpha e^{-u} \rho_{\bar{k} j})
+2n\alpha F^{j \bar{k}} R_{\bar{p} j \bar{k}}{}^{\bar q}D_{\bar q}u.
\eea
Here $R_{\bar{p} j \bar{k}}{}^{\bar q}$ denotes the curvature of the background metric $\omega$.
Introduce the notation
\bea
|DDu|_{Fg}^2=F^{j\bar k}g^{\ell\bar m}D_{j} D_{\ell} uD_{\bar k} D_{\bar{m}} u,
\quad
|D\bar Du|_{Fg}^2=F^{j\bar k}g^{\ell\bar m}D_{j} D_{\bar m}uD_{\ell}D_{ \bar k}u.
\eea
Then
\bea
2n\alpha F^{j\bar k}D_jD_{\bar k}|Du|^2
&=&
2n\alpha g^{\ell\bar m}F^{j\bar k}(D_jD_{\bar k}D_{\ell}u \,D_{\bar m}u
+
D_{\ell}u \,D_j D_{\bar k}D_{\bar m}u)
\nonumber\\
&&
+2n\alpha (|DDu|_{Fg}^2+|D\bar D u|_{Fg}^2),
\eea
and hence
\bea \label{DD_|Du|^2}
2n\alpha F^{j\bar k}D_jD_{\bar k}|Du|^2
&=&
2 \Re \langle D \sigma_2^{1/2}, Du \rangle
+2n\alpha F^{j \bar{k}} g^{l\bar p}D_l u R_{\bar pj \bar{k}}{}^{\bar q}D_{\bar q}u 
\nonumber\\
&&
- 2 \Re \langle F^{j \bar{k}} D(e^u g_{\bar{k} j} +\alpha e^{-u} \rho_{\bar{k} j}), D u \rangle 
\nonumber\\
&&
+2n\alpha (|DDu|_{Fg}^2+|D\bar D u|_{Fg}^2).
\eea
Finally, we compute the operator $2n\alpha F^{j\bar k}D_jD_{\bar k}$ acting on the Hessian $D_pD_{\bar q}u$. Differentiating the equation (\ref{D_sigma2}) again gives
\bea
F^{j\bar k}D_pD_{\bar q}g_{\bar kj}'
=
\p_p\p_{\bar q} \sigma_2^{1/2} - F^{i \bar{j},k \bar{\ell}} D_p g'_{\bar{j} i} D_{\bar{q}} g'_{\bar{\ell} k}.
\eea
Using the definition of $g'$, we obtain
\bea \label{DDDD_u}
&&2n \alpha F^{j \bar{k}}  D_j D_{\bar{k}} D_p D_{\bar{q}} u \\&=& 2 n \alpha F^{j \bar{k}} D_p D_{\bar{q}} D_j D_{\bar{k}} u  + 2 n \alpha \left(F^{j \bar{k}} R_{\bar{q} j \bar{k}}{}^{\bar{a}} u_{\bar{a} p} - F^{j \bar{k}} R_{\bar{q} p \bar{k}}{}^{\bar{a}} u_{\bar{a} j}\right) \nonumber\\
&=& F^{j \bar{k}} D_p D_{\bar{q}} g'_{\bar{k} j} - F^{j \bar{k}} D_p D_{\bar{q}} (e^u g_{\bar{k} j} + \alpha e^{-u} \rho_{\bar{k} j}) + 2 n \alpha \left(F^{j \bar{k}} R_{\bar{q} j \bar{k}}{}^{\bar{a}} u_{\bar{a} p} - F^{j \bar{k}} R_{\bar{q} p \bar{k}}{}^{\bar{a}} u_{\bar{a} j}\right) \nonumber\\
&=& \p_p\p_{\bar q} \sigma_2^{1/2} - F^{i \bar{j},k \bar{\ell}} D_p g'_{\bar{j} i} D_{\bar{q}} g'_{\bar{\ell} k} + 2 n \alpha \left(F^{j \bar{k}} R_{\bar{q} j \bar{k}}{}^{\bar{a}} u_{\bar{a} p} - F^{j \bar{k}} R_{\bar{q} p \bar{k}}{}^{\bar{a}} u_{\bar{a} j}\right) \nonumber\\
&& - F^{j \bar{k}} (e^u g_{\bar{k} j} - \alpha e^{-u} \rho_{\bar{k} j} )D_p D_{\bar{q}} u - F^{j \bar{k}} (e^u g_{\bar{k} j} + \alpha e^{-u} \rho_{\bar{k} j} )D_p u D_{\bar{q}} u \nonumber\\
&& + \alpha e^{-u} F^{j \bar{k}} D_p u D_{\bar{q}} \rho_{\bar{k} j} + \alpha e^{-u} F^{j \bar{k}} D_{\bar{q}} u D_{p} \rho_{\bar{k} j}  - \alpha e^{-u} F^{j \bar{k}} D_p D_{\bar{q}} \rho_{\bar{k} j}.
\eea

\subsection{Proof of the $C^2$ estimate}

\begin{proposition}
\label{C2estimatesec}
Let $u$ be a smooth solution to (\ref{kappa_equation}) such that $\lambda' \in \Gamma_2$, and suppose the $C^0$ estimate $B_2^{-1} M_0 \leq e^{u} \leq B_1 M_0$ holds. Suppose the parameter $\alpha<0$. There exists an $M'$ such that for all $M_0 \geq M'$, there exists $C>1$ such that
\be
\sup_X| \p \bar{\p} u |_{\o} \leq C(1 + \sup_X | Du |_{\o}^2),
\ee
where $C$ depends on $(X,\o)$, $\rho$, $\mu$, $\alpha$, $M_0$, $B_1$, $B_2$.
\end{proposition}
We will use the notation
\be
K = \sup_X |Du|^2 +1.
\ee
As before, $\lambda'=(\lambda_1', \dots, \lambda'_n)$ will denote the eigenvalues of $g'$ with respect to $g$, and we shall take the ordering $\lambda'_1 \geq \lambda'_2 \geq \cdots \geq \lambda'_n$. We will often use that the complex Hessian of $u$ can be bounded by $\lambda_1'$. Indeed, since $g' \in \Gamma_2$, we can estimate
\be \label{DDu_lambda'}
|2n \alpha u_{\bar{k} j}| \leq |g'_{\bar{k} j}| + |e^u g_{\bar{k} j} + \alpha e^{-u} \rho_{\bar{k} j}| \leq C( \lambda'_1 + 1).
\ee
\par
We now first state a lemma which exploits the specific function $w$ and the sign of the parameter $\alpha<0$.
\begin{lemma} \label{c2lemma}
Let $u$ be as Proposition \ref{C2estimatesec}. Suppose that at $p \in X$, we have $\lambda'_1 \gg K \geq 1 + |Du|^2$. Then at $p$, there holds
\be \label{large_calF}
\F \gg 1+ |Du|,
\ee
%\be
%(4 |\alpha| \kappa_c e^u)^{1/2} |Du| \leq w \leq {1 \over N} \F,
%\ee
\be
|\langle Dw, Du \rangle| \leq C \{ K \F + |DDu| |Du| \},
\ee
\be \label{DDw}
{D_1 D_{\bar{1}} w \over \lambda'_1} \geq -C \bigg\{ \F + (1+|Du|){|DDu| \over \lambda'_1} + {|D u_{\bar{1} 1}| \over \lambda'_1} \bigg\}.
\ee
\end{lemma}
{\it Proof:} We shall compute at a point where $g_{\bar{k} j} = \delta_{kj}$ and $g'_{\bar{k} j}$ is diagonal. Recall $\F = {n-1 \over 2w} \sum_i \lambda'_i$, and $\sum_i \lambda'_i = \lambda'_1 + \sigma_2^{1 \bar{1}}(\lambda') \geq \lambda'_1$. For choice of normalization $M_0 \gg 1$, by the $C^0$ estimate we have $e^{-u} \ll 1$. It follows that for $M_0 \gg 1$, $w>0$ and
\be
{1 \over C} (1 + |Du|^2) \leq w^2 \leq C (1 + |Du|^2).
\ee
Thus
\be
\F \geq {1 \over C} {\lambda'_1 \over w^2} w \geq {1 \over C}{\lambda'_1 \over K} (1 + |Du|).
\ee
Hence $\F \gg 1+ |Du|$. Next, we compute derivatives of $w^2$. 
\bea \label{Dw^2}
D_k w^2 &=&2 \kappa_c e^{2u} D_ku + 4 |\alpha|\kappa_c e^u D_k |Du|^2 + 2n \alpha^2 e^{-u}D_k \{ \tilde{\rho}^{i \bar{j}} u_i u_{\bar{j}}\} \nonumber\\
&&+ 4 |\alpha| \kappa_c e^u |Du|^2D_k u - 2n \alpha^2 e^{-u} \tilde{\rho}^{i \bar{j}} u_i u_{\bar{j}} D_k u \nonumber\\
&&+4n \alpha^2 e^{-u} \Re \langle \p u, \bar{\p} \rho \rangle_\o D_k u -4n \alpha^2 e^{-u} \Re D_k \langle \p u, \bar{\p} \rho \rangle_\o \nonumber\\
&&+ D_k \bigg\{  \alpha^2 e^{-2u} \sigma_2 (\rho) + 2n\alpha^2 e^{-u} \Delta_\o \rho +\alpha(n-1)  g^{a \bar{b}} \rho_{\bar{b} a} - {2n \alpha \over (n-2)!} \mu \bigg\}. 
\eea
Estimate
\bea
|\langle Dw,Du \rangle| &\leq& {1 \over 2w}|Du| |Dw^2| \\
&\leq& C \bigg\{ {|Du|^4 \over w} + {1+|Du| \over w} |Du| |DDu| + {1 + |Du|\over w} |Du| |D \bar{D} u| + {|Du| \over w}  \bigg\}. \nonumber
\eea
Using $C w \geq 1 + |Du|$, we obtain
\be
|\langle Dw,Du \rangle| \leq CK (1 + |Du| +{\lambda'_1 \over w}) + C |Du| |DDu| \leq C (K \F + |Du| |DDu|).
\ee
To complete the lemma, it remains to show (\ref{DDw}). Compute
\bea
D_1 D_{\bar{1}} w &=& {1 \over 2 w} \bigg\{ - {|D_1 w^2|^2 \over 2 w^2} + D_1 D_{\bar{1}} w^2 \bigg\} \nonumber\\
&=& {1 \over 2 w} \bigg\{ - {1 \over 2 w^2} \bigg| 4 \alpha \kappa_c e^u D_1 |Du|^2 \bigg|^2 - {1 \over 2 w^2}2\Re \langle 4 \alpha \kappa_c e^u  D_{\bar1} |Du|^2, R_1 \rangle \nonumber\\
&& - {|R_1|^2 \over 2 w^2} + D_1 D_{\bar{1}} w^2 \bigg\},
\eea
where
\bea
R_1 &=& 2\kappa_c e^{2u} D_1u + 2n \alpha^2 e^{-u}D_1 \{ \tilde{\rho}^{i \bar{j}} u_i u_{\bar{j}}\} + 4 |\alpha| \kappa_c e^u |Du|^2D_1 u - 2n \alpha^2 e^{-u} \tilde{\rho}^{i \bar{j}} u_i u_{\bar{j}} D_1 u \nonumber\\
&&+4n \alpha^2 e^{-u} \Re \langle \p u, \bar{\p} \rho \rangle_\o D_1 u -4n \alpha^2 e^{-u} \Re D_1 \langle \p u, \bar{\p} \rho \rangle_\o \nonumber\\
&&+ D_1 \bigg\{  \alpha^2 e^{-2u} \sigma_2 (\rho) + 2n\alpha^2 e^{-u} \Delta_\o \rho +\alpha(n-1)  g^{a \bar{b}} \rho_{\bar{b} a} - {2n \alpha \over (n-2)!} \mu \bigg\}. 
\eea
We estimate
\be
|R_1| \leq C |Du|^3 + C(n,\alpha,\rho) e^{-u}(1+|Du|) \sum_p \{ |u_{1p}| + |u_{\bar{1} p}| \} + C,
\ee
\be
|R_1|^2 \leq C |Du|^6 +  C(n,\alpha,\rho) e^{-2u} (1+|Du|^2) \sum_p \{ |u_{1p}|^2 + |u_{\bar{1} p}|^2 \} + C,
\ee
\be
\Re \langle 4 \alpha \kappa_c e^u  D_{\bar1} |Du|^2, R_1 \rangle \leq  C(n,\alpha,\rho)  (1+|Du|^2) \sum_p \{ |u_{1p}|^2 + |u_{\bar{1} p}|^2 \} + C|Du|^6 + C.
\ee
Together with $e^u \gg 1$ for $M_0 \gg 1$, we use the above inequalities to obtain
\bea
& \ & 
%|\{ 4 \alpha \kappa_c e^u  D_1 |Du|^2  \} \{ R_1 \}|
\Re \langle 4 \alpha \kappa_c e^u  D_{\bar1} |Du|^2, R_1 \rangle + {|R_1|^2 \over 2} \nonumber\\
&\leq& \bigg( {1 \over 8} (4 \alpha \kappa_c)^2 e^{2u} |Du|^2 + e^u \bigg) \sum_p \{ |u_{1p}|^2 + |u_{\bar{1} p}|^2\} + C|Du|^6 + C.
\eea
Therefore
\bea
D_1 D_{\bar{1}} w &\geq& {1 \over 2 w} \bigg\{ - {1 \over 2 w^2} (4 \alpha \kappa_c)^2 e^{2u} |Du|^2 \sum_p \bigg( {5 \over 4} |u_{1p}|^2 + 5 |u_{\bar{1} p}|^2 \bigg) \nonumber\\
&&- {1 \over w^2} \bigg( {1 \over 8} (4 \alpha \kappa_c)^2 e^{2u}|Du|^2 + e^u \bigg) \sum_p \{ |u_{1p}|^2 + |u_{\bar{1} p}|^2\} \nonumber\\
&&- {C \over w^2} |Du|^6 - C + D_1 D_{\bar{1}} w^2 \bigg\}.
\eea
Combining terms
\bea
D_1 D_{\bar{1}} w &\geq& {1 \over 2 w} \bigg\{ - {3 \over 4 w^2} \bigg( (4 \alpha \kappa_c)^2 e^{2u} |Du|^2 + e^u \bigg) \sum_p |u_{1p}|^2 \nonumber\\
&&- {C (1+|Du|^2) \over w^2} |D \bar{D} u|^2- {C \over w^2} |Du|^6 - C + D_1 D_{\bar{1}} w^2 \bigg\}.
\eea
For $e^{-u} \ll 1$, we have
\be
w^2 \geq {7 \over 8} \{ \kappa_c e^{2u} + 4 |\alpha| \kappa_c e^u  |Du|^2 \}.
\ee
Hence
\be
D_1 D_{\bar{1}} w \geq {1 \over 2 w} \bigg\{ - {6 \over 7} (4 |\alpha| \kappa_c) e^u  \sum_p |u_{1p}|^2 -C |D \bar{D} u|^2 - C |Du|^4 - C + D_1 D_{\bar{1}} w^2 \bigg\}.
\ee
Taking a second derivative of (\ref{Dw^2}), we estimate for $e^{-u} \ll 1$,
\bea
D_1 D_{\bar 1} w^2 &\geq& {6 \over 7} (4 \kappa_c |\alpha|) e^u \sum_p \{ |u_{1p}|^2 + |u_{\bar{1} p}|^2 \} \\
&& - C\{ (1+|Du|) |D u_{\bar{1} 1}| + (1+|Du|^2) (|DDu| + |D \bar{D} u|)  + |Du|^4 + 1 \}. \nonumber
\eea 
We see that the terms involving $\sum_p |u_{1p}|^2$ cancel. Hence
\be
D_1 D_{\bar{1}} w \geq - C \bigg\{ {(1+|Du|) \over w} \{ |Du_{\bar{1} 1}| + (1+|Du|) |DDu| \} + {|D \bar{D} u|^2 \over w} + {|Du|^4 +1 \over w}  \bigg\}.
\ee
Therefore, for $\lambda'_1 \gg K$,
\be
{D_1 D_{\bar{1}} w \over \lambda'_1} \geq - C \bigg\{  {|Du_{\bar{1} 1}| \over \lambda'_1} + (1+|Du|) {|DDu| \over \lambda'_1}  + {\lambda'_1 \over w} + |Du| +1  \bigg\}.
\ee
This inequality yields (\ref{DDw}). Q.E.D.
\bigskip
\par \noindent Given Lemma \ref{c2lemma}, we now prove the $C^2$ estimate. We shall use the maximum principle applied to the test function of Hou-Ma-Wu \cite{HMW}. Let $N>0$ be a large constant to be determined later. Let $L = 2n |\alpha| \sup_X |u|$. Define
\be
\psi(t) = N \log \left( 1 + {t \over 2L} \right),
\ee
for $|t| \leq L$. It follows that
\be \label{psi_identities}
{N \over L} > \psi' > {N \over 3 L}, \ \ \psi'' = - {|\psi'|^2 \over N}.
\ee
Define
\be
\phi(t) = - \log (2K - t), \ K = \sup_X |Du|^2 +1,
\ee
for $0 \leq t \leq K$. We have
\be \label{phi_bounds}
\phi'( |D u|^2) \leq {1 \over K} \leq 1, \ \ \phi'( |D u|^2) \geq {1 \over 2 K},
\ee
and the relationship
\be \label{phi_ode}
\phi'' = (\phi')^2 .
\ee
First, consider
\be
G_0(z,\xi)= \log ( g'_{\bar{j} k} \xi^k \bar{\xi}^j) - \psi(2n \alpha u) + \phi( |D u|^2 ),
\ee
for $z \in X$ and $\xi \in T_z^{1,0}(X)$ a unit vector. $G_0$ is not defined everywhere, but we may restrict to the compact set where $g'_{\bar{j} k} \xi^k \bar{\xi}^j \geq 0$ and obtain an upper semicontinuous function. Let $(p, \xi_0)$ be the maximum of $G_0$. Choose coordinates centered at $p$ such that $g_{\bar{j} k} = \delta_{jk}$ and $g'_{\bar{j} k}$ is diagonal. As before, we use the ordering $\lambda'_1 \geq \cdots \geq \lambda'_n$ for the eigenvalue of $g'$ with respect to $g$. At $p$, we have $\lambda'_1(p) = g'_{\bar{1} 1}(p)$, and $\xi_0(p) = \partial_1$. We extend $\xi_0(p)$ to a local unit vector field $\xi_0 = g_{\bar{1} 1}^{-1/2} {\partial \over \partial z^1}$.
Define the local function
\be
G (z)= \log (g_{\bar{1} 1}^{-1} g'_{\bar{1} 1}) -  \psi(2n \alpha u) + \phi( |D u|^2 ).
\ee
This function $G$ also attains a maximum at $p \in X$. We will compute at the point $p$. Covariantly differentiating $G$ gives
\be \label{c2testfunction1}
G_{\bar{j}} = { D_{\bar{j}} (e^u + \alpha e^{-u} \rho_{\bar{1} 1}) + 2n \alpha D_{\bar{j}} D_1 D_{\bar{1}} u \over {g'_{\bar{1} 1}}} + \phi' D_{\bar{j}} | D u|^2 - 2n \alpha \psi' D_{\bar{j}} u.
\ee
Covariantly differentiating $G$ a second time and contracting with $F^{i \bar{j}}$ yields
\bea \label{testfunction2}
F^{i \bar{j}} G_{\bar{j} i} &=& {2n \alpha \over \lambda'_1} F^{i \bar{j}} D_i D_{\bar{j}} D_1 D_{\bar{1}} u + {(e^u - \alpha e^{-u} \rho_{\bar{1} 1}) \over \lambda'_1} F^{i \bar{j}} D_i D_{\bar{j}} u + {(e^u + \alpha  e^{-u} \rho_{\bar{1} 1}) \over \lambda'_1} |Du|^2_F \nonumber\\
&& -{2 \alpha e^{-u} \over \lambda'_1}  \Re \{ F^{i \bar{j}} u_i (\rho_{\bar{1} 1})_{\bar{j}} \} + {\alpha e^{-u} \over \lambda'_1} F^{i \bar{j}} (\rho_{\bar{1} 1})_{\bar{j} i} - {|D g'_{\bar{1} 1}|^2_F \over \lambda'_1{}^2} + \phi' F^{i \bar{j}} D_i D_{\bar{j}} |D u|^2 \nonumber\\
&&+ \phi'' | D |D u|^2|^2_F  -2n \alpha \psi' F^{i \bar{j}} D_i D_{\bar{j}} u - (2n \alpha)^2 \psi'' |D u|^2_F.
\eea
Here we introduced the notation $| D\chi |^2_F = F^{j \bar{k}} D_j \chi D_{\bar{k}} \chi$. We first get an estimate for $F^{i \bar{j}} D_iD_{\bar{j}} D_1 D_{\bar{1}}u$ by using the identity (\ref{DDDD_u}) and noting that the complex Hessian of $u$ can be bounded by $\lambda'_1$ (\ref{DDu_lambda'}).
\be \label{A1}
2n \alpha F^{i \bar{j}} D_i D_{\bar{j}} D_1 D_{\bar{1}} u  \geq D_1 D_{\bar 1} w - F^{i \bar{j},k \bar{\ell}} D_1 g'_{\bar{j} i} D_{\bar{1}} g'_{\bar{\ell} k}  - C (1 + |Du|^2 + \lambda_1') \F.
\ee
From (\ref{F_DDu}), for suitable normalization $e^u \gg 1$ and hence $\alpha e^{-u} \rho_{\bar{k} j} \geq - g_{\bar{k} j}$, and so we have
\be \label{A2}
- 2n\alpha F^{i\bar j}D_iD_{\bar j}u \geq \F - w.
\ee
By (\ref{DD_|Du|^2}), we have
\bea \label{A3}
 F^{i \bar{j}} D_i D_{\bar{j}} |D u|^2
&\geq& {\Re \langle Dw, Du \rangle \over n \alpha}  - C (1+|Du|^2) \F + ( |D \bar{D} u|^2_{Fg} + |D D u|^2_{Fg}).
\eea
Using inequalities (\ref{A1}), (\ref{A2}), (\ref{A3}), in (\ref{testfunction2}) yields the following inequality at the maximum point $p \in X$ of $G$
\bea 
0 &\geq& {1 \over \lambda'_1} \bigg\{ D_1 D_{\bar 1} w - F^{i \bar{j},k \bar{\ell}} D_1 g'_{\bar{j} i} D_{\bar{1}} g'_{\bar{\ell} k} \bigg\} - {|D g'_{\bar{1} 1}|^2_F \over \lambda'_1{}^2 } + \phi'|D \bar{D} u|^2_{Fg} + \phi'|D D u|^2_{Fg}  \nonumber\\
&&+ \phi'' | D |D u|^2 |^2_F + {\phi' \over n \alpha} \Re \langle Dw, Du \rangle  -(2n \alpha)^2 \psi'' | D u|^2_F \nonumber\\
&& +\psi' \F -C \bigg\{ 1+ \phi' + \phi'|Du|^2 + {|Du|^2 \over \lambda'_1} \bigg\} \F - \psi' w - C.
\eea
We shall assume $\lambda'_1 \gg K \geq 1 + |Du|^2$, otherwise the estimate is complete. Since $\phi' \leq {1 \over K}$, we have
\be
1+ \phi' + \phi'|Du|^2 + {|Du|^2 \over \lambda'_1} \leq C.
\ee
Combining (\ref{large_calF}) with $w \leq C(1+|Du|)$, we obtain $w \leq {1 \over N} \F$ for $\lambda'_1$ large enough. Using Lemma \ref{c2lemma} on the terms involving derivatives of $w$, the main inequality becomes
\bea 
0 &\geq& -{1 \over \lambda'_1}  F^{i \bar{j},k \bar{\ell}} D_1 g'_{\bar{j} i} D_{\bar{1}} g'_{\bar{\ell} k} - {|D g'_{\bar{1} 1}|^2_F \over \lambda'_1{}^2} + \phi'|D \bar{D} u|^2_{Fg} + \phi'|D D u|^2_{Fg}  \nonumber\\
&& -C\phi' |DDu| |Du| - C(1+|Du|){|DDu| \over \lambda'_1} - C{|D u_{\bar{1} 1}| \over \lambda'_1} \nonumber\\
&&+ \phi'' | D |D u|^2 |^2_F  -(2n \alpha)^2 \psi'' | D u|^2_F + (\psi'-C) \F.
\eea 
Using the critical equation $DG=0$ and thus setting (\ref{c2testfunction1}) to zero, we obtain
\bea
{|D u_{\bar{1} 1}| \over \lambda'_1} &\leq& {C \over \lambda'_1}(1 +|Du|) + \phi' |D|Du|^2| + C \psi' |Du| \nonumber\\
&\leq& \F + C \phi' |Du| |DDu| + C \phi' |Du| \lambda'_1. 
\eea
In the last line we used (\ref{large_calF}) from Lemma \ref{c2lemma}. Using $\F \geq {n-1 \over 2} {\lambda'_1 \over w}$, we estimate
\be
\phi' |Du| \lambda'_1 \leq {\lambda'_1 \over \sqrt{K}} \leq C {\lambda'_1 \over w} \leq C \F.
\ee
Therefore, by using the equation $\sigma_2(\lambda') = w^2$ and (\ref{phi_bounds}),
\bea
& \ & C \bigg\{ \phi' |DDu| |Du|+ (1+|Du|){|DDu| \over \lambda'_1} + {|D u_{\bar{1} 1}| \over \lambda'_1} \bigg\} \nonumber\\
&\leq& C \bigg\{ \phi' |DDu| |Du|+ (1+|Du|){|DDu| \over \lambda'_1} +  \F \bigg\} \nonumber\\
&=& C \bigg\{ \bigg({\phi' \sigma_2^{1/2} |DDu|^2  \over \lambda'_1}\bigg)^{1/2} \, \bigg( {  \phi' |Du|^2 \lambda'_1 \over w}\bigg)^{1/2} + \bigg( {\phi' \sigma_2^{1/2} |DDu|^2 \over \lambda'_1} \bigg)^{1/2} \, \bigg( {(1+|Du|)^2 \over \phi' \lambda'_1 w}\bigg)^{1/2} +  \F \bigg\} \nonumber\\
&\leq& {\phi' \over 2} {\sigma_2^{1/2} |DDu|^2 \over n \lambda'_1} + C \bigg\{ {|Du|^2 \lambda'_1 \over K w} + {(1+|Du|)^2 K \over \lambda'_1 w} + \F \bigg\} \nonumber\\
&\leq& {\phi' \over 2} {\sigma_2^{1/2} |DDu|^2 \over n \lambda'_1} + C \F.
\eea
To obtain the last inequality, we used ${\lambda'_1 \over w} \leq {2 \over n-1}\F$, $\lambda'_1 \gg K$, $Cw \geq 1 + |Du|$, and (\ref{large_calF}). Next, we note that for any $\lambda' \in \Gamma_2$, we have the inequality $\lambda'_1 \sigma_2^{1 \bar{1}} \geq {2 \over n} \sigma_2$. This inequality is well-known, and a proof can be found for example in \cite{PPZ1}. It follows that  
\be \label{unif_ellipticity}
F^{i \bar{i}} \geq F^{1 \bar{1}} \geq {\sigma_2^{1/2} \over n \lambda'_1}.
\ee
Therefore
\be
-C \bigg\{ \phi' |DDu| |Du| + (1+|Du|){|DDu| \over \lambda'_1} + {|D u_{\bar{1} 1}| \over \lambda'_1} \bigg\} \geq -{\phi' \over 2} |DDu|^2_{Fg} - C \F.
\ee
The main inequality becomes
\bea 
0 &\geq& -{1 \over \lambda'_1}  F^{i \bar{j},k \bar{\ell}} D_1 g'_{\bar{j} i} D_{\bar{1}} g'_{\bar{\ell} k} - {|D g'_{\bar{1} 1}|^2_F \over \lambda'_1{}^2 } + \phi'|D \bar{D} u|^2_{Fg} + {\phi' \over 2} |D D u|^2_{Fg}  \nonumber\\
&&+ \phi'' | D |D u|^2 |^2_F  -(2n \alpha)^2 \psi'' | D u|^2_F + (\psi'-C) \F.
\eea 
At this point, the estimate follows from the argument of Hou-Ma-Wu \cite{HMW}. We present the argument for the sake of completeness. Since we are dealing with a $2$-Hessian equation, we will also use ideas from \cite{SX}. Before proceeding in cases, we use the critical equation $DG=0$ to notice the following estimate which holds for each fixed index $i$,
\bea \label{critical_3rd_order}
{F^{i \bar{i}} |D_i g'_{\bar{1} 1}|^2 \over \lambda'_1{}^2} &=& F^{i \bar{i}} \bigg| \phi' D_i |Du|^2 - 2 n \alpha \psi' D_i u \bigg|^2 \nonumber\\
&\leq& (1 + {1 \over 8}) (\phi')^2 F^{i \bar{i}} |D_i |Du|^2|^2 + C (\psi')^2 F^{i \bar{i}} |D_i u|^2 \nonumber\\
&\leq& (\phi')^2 F^{i \bar{i}} |D_i |Du|^2|^2 + {(\phi')^2 |Du|^2 \over 4} F^{i \bar{i}} \sum_p (|u_{i p}|^2 + |u_{\bar{i} p}|^2)  + C (\psi')^2 F^{i \bar{i}} |D_i u|^2 \nonumber\\
&\leq& \phi'' F^{i \bar{i}} |D_i |Du|^2|^2 + {\phi' \over 4} F^{i \bar{i}} \sum_p ( |u_{i p}|^2 + |u_{\bar{i} p}|^2)  + C (\psi')^2 F^{i \bar{i}} |D_i u|^2.
\eea
In the last line, we used the properties of the function $\phi$ given in (\ref{phi_bounds}) and (\ref{phi_ode}). We shall need the constants
\be \label{deltatau}
\delta= {\tau \over 4 - 3 \tau}, \ \ \tau = {1 \over 1 + N}.
\ee
Case (A): $-\lambda_n' \geq \delta \lambda_1'$. Using $F^{i \bar{j}, k \bar{\ell}} \leq 0$ by the concavity of $\sigma_2^{1/2}$ on the $\Gamma_2$ cone, $\psi''<0$, and using estimate (\ref{critical_3rd_order}) on ${|D g'_{\bar{1} 1}|^2_F \over \lambda'_1{}^2 }$, we obtain
\be
0 \geq  {\phi' \over 2} |D \bar{D} u|^2_{Fg} + {\phi' \over 4} |D D u|^2_{Fg} -C \psi'^2 | D u|^2_F + (\psi'-C) \F.
\ee
Using the assumption on the smallest eigenvalue $\lambda'_n$, we estimate for $\lambda'_1$ large enough
\be
-2n \alpha u_{\bar{n} n} = -\lambda'_n + e^u + \alpha e^{-u} \rho_{\bar{n} n} \geq \delta \lambda'_1  + e^u+ \alpha e^{-u} \rho_{\bar{n} n} \geq {\delta \lambda'_1 \over 2}. 
\ee
Hence
\be
{\phi' \over 2} |D \bar{D} u|^2_{Fg} \geq {1 \over 4K} F^{n \bar{n}} u_{\bar{n} n}^2 \geq {\delta^2 \over 16 K(2n\alpha)^2} F^{n \bar{n}} \lambda'_1{}^2.
\ee
Since $F^{n \bar{n}} \geq F^{i \bar{i}}$ for all indices $i$, we have
\be
0 \geq {\delta^2 \over 16 K(2n\alpha)^2} F^{n \bar{n}} \lambda'_1{}^2 -CK \psi'^2 F^{n \bar{n}} + (\psi'-C) \F.
\ee
By (\ref{psi_identities}), we can choose $\psi$ such that $(\psi'-C) \geq 0$. The estimate $\lambda'_1 \leq C(1+K)$ follows.
\bigskip
\par \noindent Case (B): $-\lambda_n' \leq \delta \lambda_1'$. We partition $\{1, \cdots, n \}$ into
\be
I = \{ i : F^{i \bar{i}} \leq \delta^{-1} F^{1 \bar{1}} \}, \ \ J = \{ i : F^{i \bar{i}} > \delta^{-1} F^{1 \bar{1}} \}.
\ee
Using (\ref{critical_3rd_order}) for each $i \in I$ occurring in ${|D g'_{\bar{1} 1}|^2_F \over \lambda'_1{}^2 }$, the main inequality becomes
\bea 
0 &\geq& -{1 \over \lambda'_1}  F^{i \bar{j},k \bar{\ell}} D_1 g'_{\bar{j} i} D_{\bar{1}} g'_{\bar{\ell} k} - {1 \over \lambda'_1{}^2} \sum_{i \in J} F^{i \bar{i}}  |D_i g'_{\bar{1} 1}|^2   + {\phi' \over 2} |D \bar{D} u|^2_{Fg} + {\phi' \over 4} |D D u|^2_{Fg}  \nonumber\\
&&+ \phi'' \sum_{i \in J} F^{i \bar{i}} | D_i |D u|^2 |^2  -(2n \alpha)^2 \psi'' | D u|^2_F - CK(\psi')^2 \delta^{-1} F^{1 \bar{1}} \nonumber\\
&&+ (\psi'-C) \F.
\eea 
Using (\ref{phi_ode}), $D G(p) =0$, and (\ref{psi_identities}),
\bea
\sum_{i \in J} \phi'' F^{i \bar{i}} | D_i |D u|^2 |^2 &=& \sum_{i \in J} F^{i \bar{i}} |\phi'|D u|^2_i|^2=\sum_{i \in J} F^{i \bar{i}} \left| {D_i g'_{\bar{1} 1} \over g'_{\bar{1} 1}}  - 2n \alpha \psi' u_i \right|^2 \nonumber\\
&\geq& {\tau \over \lambda'_1{}^2} \sum_{i \in J} F^{i \bar{i}} |D_i g'_{\bar{1} 1}|^2 - {\tau \over 1 - \tau} \sum_{i \in J} F^{i \bar{i}} \left|2n \alpha \psi' u_i \right|^2 \nonumber\\
&=& {\tau \over \lambda'_1{}^2} \sum_{i \in J} F^{i \bar{i}} | D_i g'_{\bar{1} 1}|^2 + {\tau N \over 1 - \tau} (2n \alpha)^2 \psi'' \sum_{i \in J} F^{i \bar{i}} | u_i|^2\nonumber\\
&=& {\tau \over \lambda'_1{}^2} \sum_{i \in J} F^{i \bar{i}} | D_i g'_{\bar{1} 1}|^2 + (2n \alpha)^2 \psi'' \sum_{i \in J} F^{i \bar{i}} | u_i|^2.
\eea 
In the last line we used the definition of $\tau$ (\ref{deltatau}). The main inequality becomes
\bea 
0 &\geq& -{1 \over \lambda'_1}  F^{i \bar{j},k \bar{\ell}} D_1 g'_{\bar{j} i} D_{\bar{1}} g'_{\bar{\ell} k} - {1-\tau \over \lambda'_1{}^2} \sum_{i \in J} F^{i \bar{i}} |D_i g'_{\bar{1} 1}|^2   + {\phi' \over 2} |D \bar{D} u|^2_{Fg} + {\phi' \over 4} |D D u|^2_{Fg}  \nonumber\\
&& - CK(\psi')^2 \delta^{-1} F^{1 \bar{1}} + (\psi'-C) \F.
\eea 
Terms involving $-\psi''>0$ were discarded.
\medskip
\par Recall that if $F(A)= f(\lambda_1, \cdots, \lambda_n)$ is a symmetric function of the eigenvalues of a Hermitian matrix $A$, then at a diagonal matrix $A$, we have (see \cite{Ball, Gerhardt}), 
\begin{eqnarray}
 \label{firstorder} F^{i \bar{j}} &=& \delta_{ij} f_i,\\
  \label{secondorder} F^{i \bar{j}, r \bar{s}} T_{i\bar j k} T_{r\bar s \bar k} &=& \sum f_{ij} T_{i\bar i k} T_{j\bar j \bar k} + \sum_{p\neq q}\frac{f_p - f_q}{\lambda_p-\lambda_q} | T_{p\bar q k}|^2,
\end{eqnarray}
where the second term on the right-hand side of (\ref{secondorder}) has to be interpreted as a limit if $\lambda_p=\lambda_q$. In our case $f(\lambda')=\sigma_2^{1/2}(\lambda')$, and we may compute
\be
f_p = {1 \over 2 \sigma_2(\lambda')^{1/2}} \sum_{k \neq p} \lambda'_k ,  \ \ \frac{f_p - f_q}{\lambda_p-\lambda_q} = - {1 \over 2 \sigma_2(\lambda')^{1/2}}.
\ee
Since $f(\lambda')=\sigma_2^{1/2}(\lambda')$ is concave, identity (\ref{secondorder}) gives us the following inequality
\be
F^{i \bar{j}, r \bar{s}} T_{i\bar j k} T_{r\bar s \bar k} \leq - {1 \over 2 \sigma_2^{1/2}} \sum_{p\neq q} | T_{p\bar q k}|^2.
\ee
We now estimate
\bea
-{1 \over \lambda'_1}  F^{i \bar{j},k \bar{\ell}} D_1 g'_{\bar{j} i} D_{\bar{1}} g'_{\bar{\ell} k} &\geq& - {1 \over 2  \lambda'_1 w} \sum_{i \neq j} | D_1 g'_{\bar{j} i}|^2 \nonumber\\
&\geq& {1 \over 2 \lambda'_1 w} \sum_{i \neq 1}  | D_i g'_{\bar{1} 1} - D_i ( e^u + \alpha e^{-u} \rho_{\bar{1} 1}) + \alpha D_1 (e^{-u} \rho_{\bar{1} i}) |^2 \nonumber\\
&\geq& {1 - {\tau \over 2} \over \lambda'_1} \sum_{i \neq 1} {1 \over 2 w} | D_i g'_{\bar{1} 1}|^2 - {C_\tau  \over \lambda'_1 w} (1+|Du|^2).
\eea
For any index $i \in J$, we have that $\lambda'_i \neq \lambda'_1$. We only keep indices $i \in J$ in the summation, and use the definitions of $J$ and case (B), to obtain
\bea
-{1 \over \lambda'_1}  F^{i \bar{j},k \bar{\ell}} D_1 g'_{\bar{j} i} D_{\bar{1}} g'_{\bar{\ell} k} &\geq& {1-{\tau \over 2} \over \lambda'_1} \sum_{i \in J}\frac{F^{i \bar{i}} - F^{1 \bar{1}} }{\lambda_1' - \lambda'_i} | D_i g'_{\bar{1} 1}|^2 - {C_\tau \over \lambda'_1 w} (1+|Du|^2) \nonumber\\
&\geq& {1 \over  \lambda'_1{}^2} (1-{\tau \over 2}) \frac{1 - \delta}{1+\delta} \sum_{i \in J} F^{i \bar{i}} | D_i g'_{\bar{1} 1}|^2 - \F \nonumber\\
&=& {1 \over  \lambda'_1{}^2} (1-\tau)  \sum_{i \in J} F^{i \bar{i}} | D_i g'_{\bar{1} 1}|^2 - \F,
\eea
for $\lambda_1' \gg K$. In the last line we used the definition of $\delta$ (\ref{deltatau}). The main inequality becomes
\be 
0 \geq {\phi' \over 2} |D \bar{D} u|^2_{Fg}  - CK(\psi')^2 \delta^{-1} F^{1 \bar{1}} + (\psi'-C) \F.
\ee 
Choosing $N \gg 1$ such that $\psi'-C \geq 0$, we obtain
\be
0 \geq {1 \over 4K} F^{1 \bar{1}} \lambda'_1{}^2 - CK (\psi')^2 \delta^{-1} F^{1 \bar{1}}.
\ee
The estimate $\lambda'_1 \leq C (1+K)$ follows. Q.E.D.

\section{The Gradient Estimate}
\setcounter{equation}{0}

The gradient estimate is immediate from the equation (\ref{kappa_equation}) when $\alpha>0$ for $e^{u} \gg 1$, since $\sigma_2(\lambda') \geq 0$. In the present case, when $\alpha<0$, we use the blow-up argument and Liouville theorem of Dinew-Kolodziej \cite{DK}.

\begin{proposition}
Let $u$ be a solution to (\ref{kappa_equation}) with parameter $\alpha < 0$ such that $\lambda' \in \Gamma_2$. Suppose the $C^0$ estimate $B_2^{-1} M_0 \leq e^{u} \leq B_1 M_0$ and $C^2$ estimate $| \p \bar{\p} u | \leq C(1 + \sup_X | Du |^2)$ hold. There exists $C>1$ such that
\be
|Du|^2 \leq C,
\ee
where $C$ depends on $(X,\o)$, $\rho$, $\mu$, $\alpha$, $M_0$, $B_1$, $B_2$.
\end{proposition}
{\it Proof:} Proceed by contradiction. Suppose there exists a sequence of functions $u_k: X \rightarrow {\bf R}$ solving (\ref{kappa_equation}) with $\lambda' \in \Gamma_2$ such that $|Du_k(x_k)| = C_k$ for some $x_k \in X$ and $C_k \rightarrow \infty$. After taking a subsequence, we may assume that $x_k \rightarrow x_\infty$ for some $x_\infty \in X$. We take a coordinate chart centered at $x_\infty$, and assume that all $x_k$ are inside this coordinate chart. We shall take coordinates such that $\o(0) = \beta$, where $\beta= \sum_j i dz^j \wedge d \bar{z}^j$. Define the local functions
\be
\hat{u}_k (x) = u_k \bigg( {x \over C_k} + x_k \bigg).
\ee
These functions have the following properties
\be
|D \hat{u}(x)| \leq |D \hat{u}_k(0)| = 1, \ \ |\p \bar{\p} \hat{u}_k| \leq {C \over C_k^2}(1 + C_k^2) \leq C, \ \ \| \hat{u}_k \|_{L^\infty} \leq C.
\ee
Elliptic estimates for the Laplacian show that $\hat{u}_k$ is bounded is $C^{1,\alpha}$. Therefore, on any $B_R(0)$ there exists a subsequence $(2n \alpha) \hat{u}_k \rightarrow u_\infty$ in $C^{1,\beta}(B_R(0))$. Define $\Phi_k: \C^n \rightarrow \C^n$ to be the map $\Phi_k(x) = C_k^{-1} x + x_k$. We introduce notation analogous to \cite{TW},
\be
\beta_k = C_k^2 \, \Phi_k^* \o, \ \ \chi_k = \Phi_k^* (e^{u_k} \o + \alpha e^{-u_k} \rho).
\ee
It follows that
\be
\Phi_k^*(e^{u_k} \o + \alpha e^{-u_k} \rho + 2 n \alpha i \ddb u_k) = \chi_k + 2n \alpha \, i \ddb \hat{u}_k.
\ee
We also note the following convergence
\be \label{beta_chi}
\beta_k \rightarrow \beta, \ \ \chi_k \rightarrow 0, \ {\rm in} \ C^\infty_{{\rm loc}}.
\ee
Since $\chi_k + 2n \alpha i \ddb \hat{u}_k$ is in the $\Gamma_2$ cone, it follows that for any function $v \in C^2(B_R(0))$ such that $(i \ddb v) \in \Gamma_2$, we have $ (\chi_k + 2n \alpha i \ddb \hat{u}_k) \wedge i \ddb v \wedge \beta^{n-2} \geq 0$. This follows from Garding's inequality $\sum_i {\p \sigma_2 (\lambda) \over \p \lambda_i} \mu_i \geq 2 \sigma_2(\lambda)^{1/2} \sigma_2(\mu)^{1/2}$ for any $\lambda, \mu \in \Gamma_2$. Hence upon taking a limit, we have
\be
i \ddb u_\infty \wedge i \ddb v \wedge \beta^{n-2} \geq 0
\ee
in the sense of currents. This is the definition of $2$-subharmonicity introduced by Blocki \cite{Bl} that is required in the Liouville theorem of Dinew-Kolodziej \cite{DK}. Having shown that $u_\infty$ is $2$-subharmonic, we will show that it is maximal by proving $(i \ddb u_\infty)^2 \wedge \beta^{n-2} =0$ as a measure. After multiplying through by $(C_k^2)^{n-2}$ and pulling back by $\Phi_k$, the equation (\ref{kappa_equation}) solved by $u_k$ becomes the following equation for $\hat{u}_k$
\bea
& \ & (\chi_k + 2n\alpha i \ddb \hat{u}_k)^2 \wedge \beta_k^{n-2} \\
&=& \bigg\{ {\kappa_c e^{2\hat{u}_k} \over C_k^4}   -4 \alpha \kappa_ce^{\hat{u}_k} {|D\hat{u}_k|^2 \over C_k^2} + 2n \alpha^2 e^{-\hat{u}_k} {\tilde{\rho}^{a \bar{b}} D_a (\hat{u}_k) D_{\bar{b}} (\hat{u}_k) \over C_k^2} - {4n \alpha^2 \over C_k^3} e^{-\hat{u}_k} \Re \langle \p \hat{u}_k, \bar{\p} \rho \rangle_\o  \nonumber\\
&& + {\alpha^2 \over C_k^4} e^{-2\hat{u}_k} \sigma_2 (\rho) + {2n\alpha^2 \over C_k^4} e^{-\hat{u}_k} \Delta_\o \rho + {\alpha(n-1) \over C_k^4}  g^{a \bar{b}} \rho_{\bar{b} a} - {2n \alpha \over C_k^4 } {\mu \over (n-2)!} \bigg\} (\o \circ \Phi_k)^2 \wedge \beta_k^{n-2}.\nonumber
\eea
Since $\hat{u}_k$ is uniformly bounded, $|D \hat{u}_k| \leq 1$ and $C_k \rightarrow \infty$, we see that the right hand side tends to zero. Combining this with (\ref{beta_chi}), we may conclude
\be
(i \ddb \hat{u}_k)^2 \wedge \beta^{n-2} \rightarrow 0,
\ee
in the sense of currents. Since $2n \alpha \, \hat{u}_k \rightarrow u_\infty$ locally uniformly, it is well-known (e.g. \cite{De} Chapter III Cor. 3.6) that $(2n \alpha \, i \ddb \hat{u}_k)^2 \wedge \beta^{n-2} \rightarrow (i \ddb u_\infty)^2 \wedge \beta^{n-2}$ weakly. Thus
\be
( i \ddb u_\infty)^2 \wedge \beta^{n-2}= 0,
\ee
in the sense of Bedford-Taylor \cite{BT1}. Since $u_\infty$ is a bounded maximal $2$-subharmonic function in $\C^n$ with bounded gradient, by the Liouville theorem of Dinew-Kolodziej \cite{DK}, $u_\infty$ must be constant. We obtain a contradiction, since $|D u_\infty|^2(0) = 1$. Q.E.D.

\section{Solving the Fu-Yau equation}
\setcounter{equation}{0}
We return to the continuity method (\ref{continuity_meth1})
\be \label{continuity_meth2}
i \ddb (e^{u_t} \o - t \alpha e^{-u_t} \rho) \wedge \o^{n-2} + n \alpha i \ddb u_t \wedge i \ddb u_t \wedge \o^{n-2} + t \mu {\o^n \over n!} =0. 
\ee
We combine our estimates to establish
\begin{proposition} \label{apriori_est}
Let $\alpha<0$. There exists $M' \gg 1$ such that for all $M_0 \geq M'$, the following holds. Let $u_0 = \log M_0$, and suppose that for all $t \in [0,t_0)$ with $t_0 \leq 1$ there exists a solution $u_t$ to (\ref{continuity_meth2}) such that $\lambda'_{(t,u_t)} \in \Gamma_2$ and $\int_X e^{u_t}=M_0$. Then there exist constants $C>1$ and $0<\gamma<1$ only depending on $(X,\o)$, $\rho$, $\alpha$, $\mu$ and $M_0$ such that
\be
\| u_t \|_{C^{2,\gamma}} \leq C, \ \ \sigma_2(\lambda'_{(t,u_t)}) \geq {1 \over C}.
\ee
%The constants $M_0$ and $C$ only depend on $(X,\o)$, $\rho$, $\alpha$ and $\mu$.
\end{proposition}
{\it Proof:} Combining our $C^0$, $C^1$ and $C^2$ estimates yields $\| u_t \|_{C^2} \leq C$. When $\alpha<0$, from (\ref{kappa_equation}) it is clear that for $M_0 \gg 1$ we have a lower bound for $\sigma_2(\lambda')$. From (\ref{unif_ellipticity}), we see that $\sigma_2(\lambda')^{1/2}$ is a concave uniformly elliptic operator, with right-hand side in $C^{\gamma'}$ for some $0<\gamma'<1$. We may apply the same argument as in \cite{PPZ1} using a result by Tosatti-Wang-Weinkove-Yang \cite{TWWY}, extending an original argument of Wang \cite{W}, to obtain $\| u_t \|_{C^{2,\gamma}} \leq C$. Q.E.D.
\bigskip
\par \noindent Since $I$ is open and contains $0$, it contains the interval $[0, \hat{t})$. Consider any sequence $t_k \in I$ converging to $\hat{t}$. Then there exists $(t_k, u_k) \in B_1$ satisfying (\ref{continuity_meth1}), and by Proposition \ref{apriori_est}, the estimate $\| u_k \|_{C^{2,\gamma}} \leq C$ holds. By Arzela-Ascoli, after passing to a subsequence we have that $u_k \rightarrow \hat{u}$ in $C^{2,\gamma}$. By the lower bound for $\sigma_2(\lambda')$ in Proposition \ref{apriori_est}, we have $\lambda'_{(\hat{t}, \hat{u})} \in \Gamma_2$. Taking the limit yields $(\hat{t}, \hat{u}) \in B_1$ and $\Psi(\hat{t}, \hat{u})=0$, hence $I$ contains $\hat{t}$. It follows that $I$ is closed, and together with the fact that it is already known to be open and not empty, that $I = [0,1]$.

\medskip 
We have shown the existence of a $C^{2,\gamma}(X,\R)$ solution to (\ref{continuity_meth1}) at $t=1$. Differentiating $\sigma_2(\lambda')^{1/2}$ yields (\ref{LDu1}) and (\ref{LDu2}). We know that $F^{j \bar{k}} D_j D_{\bar{k}}$ is uniformly elliptic with coefficients in $C^\gamma$. Since $\sigma_2(\lambda')$ is a smooth function of $(z,u,Du)$, we have that $\p_p \sigma_2^{1/2}$ is in $C^\gamma$. By Schauder estimates and a bootstrapping argument, we see that this solution $u$ is smooth. This establishes existence of solutions to the Fu-Yau equation when $\alpha<0$. 

\section{Application to the Strominger system} \label{strominger_syst_section}
\setcounter{equation}{0}

By a construction of Fu and Yau \cite{FY1, FY2}, solutions of the Fu-Yau equation can be viewed as particular solutions of the Strominger system. Our goal is this section is to describe briefly a specific example due to \cite{FY2}, which satisfies $\alpha=-2$ and $\mu=0$, so that Theorem \ref{main_theorem} is directly applicable.

\medskip
Let $X$ be a compact Calabi-Yau manifold of dimension $n$ with Ricci-flat metric $\o_X$ and nowhere vanishing holomorphic $(n,0)$ form $\Omega_X$. Let ${\o_1 \over 2 \pi}, {\o_2 \over 2 \pi} \in H^2(X, \Z)$ be primitive harmonic $(1,1)$ forms. This data determines a $T^2$ fibration $\pi: M \rightarrow X$ with a $1$-form $\theta$, such that for every $u \in C^\infty(X, \R)$, the Hermitian form 
\be
\o_u = \pi^*(e^u \o_X) + {i \over 2} \theta \wedge \bar{\theta}, 
\ee
is a metric $\o_u>0$ on $M$, and
\be
\Omega = \Omega_X \wedge \theta,
\ee
is a nowhere vanishing holomorphic $(n+1,0)$ form. Furthermore, the balanced condition (\ref{balanced}) is satisfied for any $\o_u$. If we take a stable vector bundle $E$ over $(X,\o_X)$ with Hermitian-Einstein metric $H$ (whose existence is guaranteed by the Donaldson-Uhlenbeck-Yau theorem), then the system $(\pi^*(E), \pi^*(H),M,\o_u)$ automatically satisfy the conditions (\ref{HE}) and (\ref{balanced}) in the Strominger system. Thus it remains only to solve the last condition (\ref{anomaly}), which is equivalent, as explained earlier in the Introduction, to the Fu-Yau equation (\ref{FY_eqn}) on $X$ with $\mu$ given explicitly by (\ref{mu}). (The exact form of the form $\rho$ is not needed for our considerations.)

\medskip

The following explicit example with $\alpha=-2$ and $\mu=0$ was given in \cite{FY2}.  
Choose line bundles $L_1$, $L_2$ equipped with metrics $h_1$, $h_2$, such that their curvature forms satisfy $i F_{h_1} = -i \ddb \log h_1 = \o_1$ and $i F_{h_2} = - i \ddb \log h_2 = \o_2$. Let
\be
E = L_1 \oplus L_2 \oplus T^{(1,0)}X, \ \ H = (h_1,h_2, \o_X).
\ee
The curvature of $(E,H)$ is $i F_H = {\rm diag}(\o_1, \o_2, i R_X)$. Since $\o_1$ and $\o_2$ are primitive and $R_X$ is Ricci-flat, we have $F_H \wedge \o_X^{n-1} =0$. Using that for a primitive $(1,1)$ form $\eta$ we have $*_{\o_X} \eta = - {1 \over (n-2)!} \o_X^{n-2} \wedge \eta$, we compute
\bea
{\alpha \over 4} {\rm Tr}( F_H \wedge F_H - R_X \wedge R_X) \wedge \o^{n-2}_X &=& -{1 \over 2} ( - \o_1 \wedge \o_1 - \o_2 \wedge \o_2) \wedge \o_X^{n-2} \nonumber\\
&=& -{(n-2)! \over 2} (\|\o_1 \|^2_{\o_X} + \| \o_2 \|^2_{\o_X}) {\o_X^n \over n!}.
\eea
It follows that $\mu=0$ in this case. The system $(\pi^*E, \pi^*H, M, \o_u)$ satisfies the modified Strominger system.

\section{Other generalizations of the Strominger system
in terms of higher Chern classes}
\setcounter{equation}{0}

Finally, we would like to observe that another natural generalization of the Strominger system may be
\be \label{HE_bc}
F_H \wedge \o^{n} = 0, \ \ F_H^{2,0} = F_H^{0,2} = 0
\ee
\be
\label{anomaly_bc}
i \ddb(\|\Omega\|_{\o}^{2(n-2)\over n} \o^{n-1}) - {\alpha \over 4} ( {\rm Tr}\, R \wedge\cdots\wedge R - {\rm Tr}\, F_H \wedge\cdots\wedge F_H) = 0,
\ee
\be \label{balanced_bc}
d \bigg( \| \Omega \|_{\o}^{{2(n-1) \over n}} \o^{n} \bigg) = 0.
\ee
Here $E\to M$ is a holomorphic vector bundle over a compact complex manifold of dimension $n+1$, and $H$ and $\o$ are Hermitian metrics on $E$ and on $M$ respectively. The left-hand side in (\ref{anomaly_bc}) is an $(n,n)$-form, so that the curvatures $R$ and $F$ in each wedge product appears $n$ times, giving the $n$-th Chern classes for $T^{1,0}(M)$ and for $E$ respectively.
On a Goldstein-Prokushkin fibration, it is easy to see that the Hermitian metrics $\o_u$ satisfy the conformally balanced condition (\ref{balanced_bc}) for any smooth function $u$ on the $n$-dimensional base Calabi-Yau manifold $X$. By choosing as before a stable bundle $E$ with its corresponding Hermitian-Einstein metric $H$, we can then reduce this system to the sole equation (\ref{anomaly_bc}). This equation is in turn a scalar equation involving complex Hessian operators which may be of interest in itself. The above system also leads to a natural generalization to arbitrary dimension $n+1$ of the anomaly flow defined in \cite{PPZ3} for $n+1=3$. We shall return to these issues elsewhere.

\bigskip
\noindent
{\bf Acknowledgements}: The authors are very grateful to Li-Sheng Tseng for
discussions on modifications of the Strominger system.

\bigskip
Department of Mathematics, Columbia University, New York, NY 10027, USA

\smallskip

phong@math.columbia.edu

\bigskip
Department of Mathematics, Columbia University, New York, NY 10027, USA

\smallskip
picard@math.columbia.edu

\bigskip
Department of Mathematics, University of California, Irvine, CA 92697, USA

\smallskip
xiangwen@math.uci.edu

\end{document}